\documentclass[preprint]{elsarticle}

\usepackage{amsmath}
\usepackage{amsthm}
\usepackage{amssymb}		
\usepackage{amsfonts}
\usepackage{eucal}		

\usepackage{arydshln}           
\renewcommand\ge\geqslant
\renewcommand\le\leqslant

\newtheorem{theorem}{Theorem}[section]

\newtheorem{lemma}[theorem]{Lemma}

\theoremstyle{definition}

\theoremstyle{remark}

\numberwithin{equation}{section}

\numberwithin{equation}{section}

\newcommand\myatop[2]{\genfrac{}{}{0pt}{}{#1}{#2}}

\newcommand{\defeq}{\mathrel{:=}}

\makeatletter
\newenvironment{sizeequation}[1]{%
  \skip@=\baselineskip
  #1%
  \baselineskip=\skip@
  \equation
}{\endequation \ignorespacesafterend}
\newenvironment{sizedisplaymath}[1]{%
  \skip@=\baselineskip
  #1%
  \baselineskip=\skip@
  \displaymath
}{\enddisplaymath \ignorespacesafterend}

\newenvironment{sizealign}[1]{%
  \skip@=\baselineskip
  #1%
  \baselineskip=\skip@
  \align
}{\endalign \ignorespacesafterend}  
\newenvironment{sizemultline}[1]{%
  \skip@=\baselineskip
  #1%
  \baselineskip=\skip@
  \multline
}{\endmultline \ignorespacesafterend}  
  
\makeatother

\begin{document}

\begin{frontmatter}
\title{Extensions of the Classical Transformations of the Hypergeometric Function~${}_3F_2$}
\author{Robert S. Maier}
\ead{rsm@math.arizona.edu}
\address{Depts.\ of Mathematics and Physics, University of Arizona, Tucson,
AZ 85721, USA}
\begin{abstract}
It is shown that the classical quadratic and cubic transformation
identities satisfied by the hypergeometric function ${}_3F_2$ can be
extended to include additional parameter pairs, which differ by
integers.  In the extended identities, which involve hypergeometric
functions of arbitrarily high order, the added parameters are
nonlinearly constrained: in the quadratic case, they are the negated
roots of certain orthogonal polynomials of a discrete argument (dual
Hahn and Racah ones).  Specializations and applications of the
extended identities are given, including an extension of Whipple's
identity relating very well poised ${}_7F_6(1)$ series and balanced
${}_4F_3(1)$ series, and extensions of other summation identities.
\end{abstract}
\end{frontmatter}

\section{Introduction}
\label{sec:intro}

The Gauss hypergeometric function~${}_2F_1$ and its non-confluent
generalizations of higher order, such as ${}_3F_2$, ${}_4F_3$, etc.,
are parametric higher transcendental functions of continuing
importance in pure and applied mathematics.  As a function of a
complex argument~$x$, each is defined as a parametric series that
converges on the unit disk centered on $x=0$.  These functions satisfy
many transformation identities of the form
$F(\varphi(x))=\mathcal{A}(x)\tilde F(x)$, where $\varphi$~is a
rational function satisfying $\varphi(0)=0$, $\mathcal{A}$~is a
product of zero or more powers of rational functions, and the
parameters of the left-hand hypergeometric function~$F$ and its lifted
version~$\tilde F$ are constrained and related.  The best known
identities of this type are Euler's and Pfaff's transformations
of~${}_2F_1$, for which $\varphi$~is of degree~$1$, and the quadratic
and cubic transformations of~${}_2F_1$.  A~longer list of
transformations of~${}_2F_1$ with at~least one free parameter was
obtained by Goursat~\cite{Goursat1881}.

Only a few of the transformations of~${}_2F_1$ to itself extend to
ones of ${}_3F_2$ to itself~\cite{Askey94}.  On the ${}_3F_2$ level,
the classical identities include Whipple's quadratic
transformation~\cite[(3.1.15)]{Andrews99} and Bailey's two cubic
ones~\cite[Ch.~3, Ex.~3.8]{Andrews99}.  In each, the left-hand
${}_3F_2$ has parametric excess equal to~$\frac12$.  (The parametric
excess or Saalsch\"utzian index of a hypergeometric function is the
sum of its lower parameters, less the sum of its upper ones;
throughout this paper, it will be denoted by~$S$.)  Each of these
three has a `companion' in which the left-hand function~$F$ has
$S=-\frac12$ and the right-hand function~$\tilde F$ is not a~${}_3F_2$
but a~${}_4F_3$.  (See \cite[p.~97,\ Example~6]{Bailey35} and
\cite[(4.1),(5.4),(5.7)]{Gessel82}.)

If the hypergeometric functions $F,\tilde F$ are of like order, a
transformation of the form $F(\varphi(x))=\mathcal{A}(x)\tilde F(x)$
may be attributable to the differential equation satisfied by~$F$
being lifted by (i.e., pulled back along) the map $x\mapsto
\varphi(x)$, to the differential equation satisfied by~$\tilde F$.
(For the case of~${}_2F_1$, see \cite[\S\,3.9]{Andrews99}
and~\cite{Vidunas2009}.)  Recently, Kato determined all
transformations of ${}_3F_2$ to~${}_3F_2$ which are of this
type~\cite{Kato2008}.  They include Whipple's quadratic, Bailey's two
cubics, and several more obscure ones.

In this paper, hypergeometric function transformations of a novel kind
are investigated: identities of the form
$F(\varphi(x))=\mathcal{A}(x)\tilde F(x)$ in which the order
of~$\tilde F$ may be arbitrarily larger than the order of~$F$.  In
particular, it is shown that each of the three classical
transformations of a ${}_3F_2$ (with $S=\frac12$) to another~${}_3F_2$
can be extended to one of a ${}_3F_2$ (with $S=\frac12+k$,
$k=0,1,2,\dots$) to a~${}_{3+2k}F_{2+2k}$.  The parameters of the
latter function,~$\tilde F$, are \emph{nonlinearly constrained}: they
arise from the (negated) roots of a certain polynomial.  An example is
the extension of Whipple's quadratic, which is
{
  \setlength\arraycolsep{2.0pt}
  \begin{equation}
    \label{eq:Niblettplus}
    \begin{split}
      &{}_3F_2\left[
        \begin{array}{c}
          {\tfrac{a}2,\:\:\tfrac12+\tfrac{a}2,\:\:1-k+a-b-c} \\
          {1+a-b,\:\: 1+a-c}
        \end{array}
        \biggm|
        -\,\frac{4x}{(1-x)^2}
        \right] \\
      &\qquad=(1-x)^a\: {}_{3+2k}F_{2+2k}
      \left[
        \begin{array}{ccc:ccc}
          a, & b, & c, & 1+\xi_1, & \dots, & 1+\xi_{2k}\\
          & 1+a-b, & 1+a-c, & \xi_1, & \dots, & \xi_{2k}
        \end{array}
        \biggm| x
        \,
        \right].
    \end{split}
  \end{equation}
}The notable feature of the right-hand ${}_{3+2k}F_{2+2k}$ is that it has $2k$ unit-difference parameter-pairs, the lower parameters $\xi_1,\dots,\xi_{2k}$ of which are constrained
to equal the negated roots of
\begin{equation}
  Q_k^{(2)}(n;\,a;\,b,c) = {}_3F_2
  \left[
    \begin{array}{c}
      {-n,\:n+a,\:-k} \\
      {b,\:c}
    \end{array}
    \biggm|1\,
    \right],
\end{equation}
which is a polynomial of degree~$2k$ in~$n$.  

This result makes unexpected contact with the theory of orthogonal
polynomials of a discrete variable, because $Q_k^{(2)}(n;\,a;\,b,c)$
is essentially a dual Hahn polynomial~\cite[\S\,9.6]{Koekoek2010}: it
is invariant under $n\mapsto -n-a$ and can be written as
$R_k(\lambda(n);a;b,c)$, where $R_k(\lambda;a;b,c)$ is of degree~$k$
in $\lambda(n)=n(n+a)$, the so-called coordinate of a quadratic
lattice.  The case $k=0$ of~(\ref{eq:Niblettplus}) is the classical
one; the case $k=1$ was proved more recently~\cite{Niblett51}, as was
its $q$\nobreakdash-analogue~\cite{AlSalam84}.  It should be noted
that for all $k\ge0$, the ${}_{3+2k}F_{2+2k}$
in~(\ref{eq:Niblettplus}), having $2k$ unit-difference parameter-pairs, can be written as a finite sum of certain
${}_3F_2$'s~\cite{Karlsson71}.  But this fact is not used in the
proof.

The two cubic transformations of Bailey can be extended to $k\ge0$ in
the same way, though the corresponding degree\nobreakdash-$2k$
polynomials $Q_k^{(3)}, Q_k^{(3')}$ are asymmetric and may lack an
interpretation as orthogonal polynomials.  One of the resulting
identities is the curious specialization
\begin{sizeequation}{\small}
\begin{split}
  &{}_3F_2\left[
    {
      \myatop
      {-\tfrac16 + \tfrac{\sqrt3}3\sin\theta,
        \:\:
        \tfrac16 + \tfrac{\sqrt3}3\sin\theta,      
        \:\:
        \tfrac12 + \tfrac{\sqrt3}3\sin\theta}
      {1+\sin(\theta+\tfrac{\pi}6), \:\: 1+\sin(\theta-\tfrac{\pi}6)}
    }
    \biggm|
    -\,\frac{27x}{(1-4x)^3}
    \right]
  =(1-4x)^{-\frac12+\sqrt3\sin\theta}
  \\
  &\quad\qquad{}\times
  {}_4F_3\left[
    \begin{array}{ccc:c}
      -\tfrac12+\sqrt{3}\sin\theta, & -\tfrac12-\cos\theta, & -\tfrac12 + \cos\theta, & \tfrac32 + \tfrac{\sqrt{3}}3 \sin\theta \\
      & 1+\sin(\theta+\tfrac{\pi}6), & 1+\sin(\theta-\tfrac{\pi}6), & -\tfrac12 + \tfrac{\sqrt{3}}3 \sin\theta
    \end{array}
    \biggm|
    x\,
    \right].
\end{split}
\label{eq:curious}
\end{sizeequation}
The left-hand ${}_3F_2$ has $S=\tfrac32=\tfrac12+k$ with $k=1$.  One
would expect the right-hand function to be ${}_{3+2k}F_{2+2k} =
{}_5F_4$, but the left-hand parameters are chosen here in such a way
that the right-hand pairs
$\left[\begin{smallmatrix}1+\xi_1, & 1+\xi_2\\ \xi_1, & \xi_2\end{smallmatrix}\right]$    that come from the negated roots $\xi_1,\xi_2$ of~$Q^{(3)}_1$
    satisfy $\xi_2=1+\xi_1$.  This makes possible their merging into
    the single final pair seen in~(\ref{eq:curious}), which is of the
    form
    $\left[\begin{smallmatrix}2+\xi_1\\ \xi_1\end{smallmatrix}\right]$.

It is also shown that
the identities extending Whipple's quadratic transformation and
Bailey's cubic ones have generalizations to ${}_4F_3$.
In each identity a new parameter-pair $\left[\begin{smallmatrix}k+d
  \\ d\end{smallmatrix}\right]$, with $d$ supplying a degree of
  freedom, can be added to the parameter array of the
  left-hand~${}_3F_2$, converting it to a~${}_4F_3$.  The resulting
  generalized polynomials $Q_k^{(2)}, Q_k^{(3)}, Q_k^{(3')}$ on the
  right-hand side depend on~$d$ and have representations in~terms
  of~${}_4F_3$, and the latter two are now of degree~$3k$ in~$n$.  The
  generalized $Q_k^{(2)}$ is essentially a Racah
  polynomial~\cite[\S\,9.2]{Koekoek2010}.  The results of this paper
  on quadratic ${}_3F_2$ and ${}_4F_3$ transformations make contact
  with work of Miller and Paris~\cite{Miller2013} and Rathie, Rakha
  \emph{et al.}~\cite{Rakha2011,Rakha2009,Wang2013}, who have
  considered the effects of adding some number $r\ge1$ of
  parameter-pairs with integral differences, such as
  $\left[\begin{smallmatrix} m_1+d_1, & \dots & m_r+d_r \\ d_1, & \dots
    & d_r \end{smallmatrix}\right]$, to the left-hand functions in
  quadratic transformations of~${}_2F_1$.

On the ${}_3F_2$ level, it is additionally shown that each of the
`companion' transformations of a parametric ${}_3F_2$ with
$S=-\frac12$ to a~${}_4F_3$ (i.e., the companions of Whipple's
quadratic and Bailey's cubics) has an extension from $k=0$ to $k\ge0$.
Each extends to a transformation of a parametric ${}_3F_2$ with
$S=-\frac12-k$, $k=0,1,2,\dots$, to a ${}_{4+4k}F_{3+4k}$.  The
parameter arrays of the latter function~$\tilde F$ include $1+\nobreak
4k$ parameter-pairs with unit differences, of the form
$\left[\begin{smallmatrix} 1+\xi_1, & \dots & 1+\xi_{1+4k} \\ \xi_1, &
  \dots & \xi_{1+4k}\end{smallmatrix}\right]$.  Here,
$\xi_1,\dots,\xi_{1+4k}$ are the negated roots of a new polynomial
$\mathbf{Q}_k^{(2)}$, resp.\ $\mathbf{Q}_k^{(3)}$,
resp.~$\mathbf{Q}_k^{(3')}$.  These $k$\nobreakdash-indexed
polynomials have no obvious hypergeometric representation or
interpretation involving orthogonality, but recurrences for them are
given.  Interestingly, the new family $\mathbf{Q}_k^{(2)}$, like the
dual Hahn and Racah ones denoted by~${Q}_k^{(2)}$, is defined on a
quadratic lattice.

Gessel and Stanton~\cite{Gessel82} showed that by pairing ${}_3F_2$
transformations with their companions, one can derive many
hypergeometric evaluation formulas, including Whipple's summation
identity relating very well poised ${}_7F_6(1)$ series and balanced
${}_4F_3(1)$ series, and `strange' evaluations discovered by Gosper.
Applying the same technique to the extensions of this paper yields
extended versions of several of the Gessel--Stanton formulas, which
incorporate parameter-pairs with integral differences.  These new
formulas, in particular two extensions of Whipple's identity with
extension parameter $k=0,1,2,\dots$, overlap those recently found by
Srivastava, Vyas and Fatawat~\cite{Srivastava2018}.

Finally, a classical technique (multiplying both sides of a
hypergeometric transformation formula by a power of $(1-\nobreak x)$
and equating the coefficients of~$x^m$ on the two sides), applied to
the extensions of this paper, is shown to yield extensions of certain
summation identities due to Bailey \cite[\S\,4.5(1,2)]{Bailey35}.
Again, these are extensions from a classical case ($k=0$) to
$k=0,1,2,\dots$.


The main extension theorems are stated in \S\,\ref{sec:maintheorems},
and most are proved in~\S\,\ref{sec:proofs}.  The recurrences
satisfied by the $Q_k$ and~$\mathbf{Q}_k$, which resemble and include
those satisfied by the dual Hahn and Racah polynomials, are derived
in~\S\,\ref{sec:polys}.  The summation identities mentioned in the two
preceding paragraphs are derived in
\S\S~\ref{sec:summations1}~and~\ref{sec:summations2}.

\section{Preliminaries}
\label{sec:prelims}

The generalized hypergeometric function $F={}_{r+1}F_r$, with
$(a)=a_0,a_1,\dots,a_r$ and $(b)=b_1,\dots,b_r$ as its arrays of
$\mathbb{C}$-valued parameters, is defined by
\begin{equation}
\label{eq:hgdef}
  F
  \left[
  \begin{array}{cccc}
      a_0, & a_1, & \dots, & a_r \\
      & b_1, & \dots, & b_r 
  \end{array}
  \biggm|
  x
  \,
  \right]
  =\sum_{n=0}^\infty
  \frac{(a_0)_n(a_1)_n\dots(a_r)_n}{(1)_n(b_1)_n\dots(b_r)_n} \,x^n,
\end{equation}
the Pochhammer symbol $(c)_n$ denoting $(c)(c+1)\dots(c+n-1)$, with
$(c)_0=1$.  It is assumed that no lower parameter is a nonpositive
integer, to avoid division by zero; and if an upper one is a
nonpositive integer, the series will terminate.  The series converges
on $|x|<1$, and at $x=1$ if ${\textrm{Re}}\,S>0$; if $x=1$, the
argument is usually omitted.  Hypergeometric identities of the form
$F(\varphi(x))=\mathcal{A}(x)\tilde F(x)$ with $\varphi(0)=0$ are
taken to hold on the largest neighborhood of $x=0$ to which both sides
can be analytically continued.

Any ${}_{r+1}F_r$ with parametric excess~$S$ is said to be
$S$\nobreakdash-balanced.  It is called well-poised if $a_0+\nobreak1
= a_1+\nobreak b_1=\dots=a_r+\nobreak b_r$, or if the same holds when
$a_0,\dots,a_r$ and $b_1,\dots,b_r$ are suitably permuted, and nearly
poised if a single one of these $r+1$ parameter-pair sums differs from
the others.  It is called $(M,N)$-poised if $Ma_0+\nobreak N =
Ma_1+\nobreak Nb_1=\dots=Ma_r+\nobreak Nb_r$, where $M,N$ are positive
integers.  It is called very well poised if it is well-poised and a
parameter-pair, e.g.,
$\left[\begin{smallmatrix}a_1\\b_1\end{smallmatrix}\right]$, equals
$\left[\begin{smallmatrix}1+\frac{a_0}2\\\frac{a_0}2\end{smallmatrix}\right]$.

    It is convenient to extend the definition~(\ref{eq:hgdef}) to
\begin{equation}
\label{eq:hgdef2}
  F
  \left[
  \begin{array}{c:c|c}
    {\myatop{\textstyle(\alpha)}{\textstyle(\beta)}} & \, Q(n) & \, x
  \end{array}
  \right]
  =\sum_{n=0}^\infty
  \frac{((\alpha))_n}{(1)_n((\beta))_n}\,Q(n) \,x^n,
\end{equation}
where $(\alpha),(\beta)$ are arrays of parameters, with $((\alpha))_n
\defeq \prod_i(\alpha_i)_n$ as usual, and
$Q\colon\mathbb{N}\to\mathbb{C}$ is any weighting function of growth
no more rapid than exponential.  If $Q(n)$ is a polynomial of
degree~$\ell$ satisfying $Q(0)=1$, with $(\alpha)=a_0,a_1,\dots,a_r$
and $(\beta)=b_1,\dots,b_r$, it follows from (\ref{eq:hgdef})
and~(\ref{eq:hgdef2}) that
\begin{equation}
\label{eq:hypg}
\begin{split}
  &F
  \left[
    \begin{array}{c:c|c}
      \begin{array}{cccc}
        a_0, & a_1, & \dots, & a_r \\
        & b_1, & \dots, & b_r 
      \end{array}
      &
      \,
      Q(n)
      &
      \,
      x
    \end{array}
    \right]
  \\
  &\qquad={}_{r+\ell+1}F_{r+\ell}
  \left[
  \begin{array}{cccc:ccc}
      a_0, & a_1, & \dots, & a_r, & 1+\xi_1, & \dots, & 1+\xi_\ell \\
      & b_1, & \dots, & b_r, & \xi_1, & \dots, & \xi_\ell
  \end{array}
  \biggm|
  x
  \,
  \right],
\end{split}
\end{equation}
where $\xi_1,\dots,\xi_\ell$ are the negated roots (i.e., zeroes)
of~$Q(n)$, counted with multiplicity.  The right-hand side
of~(\ref{eq:hypg}) is a hypergeometric function with
$\ell$~unit-difference parameter-pairs.  It can be obtained from the
function of~(\ref{eq:hgdef}) by acting on~it with the differential
operator $Q(\delta)$, where $\delta = x\frac{{\rm d}}{{\rm d}x}$.  In
the formulas that employ the notation of~(\ref{eq:hgdef2}), the
normalization $Q(0)=1$ will hold, with one exception to be noted.

Any hypergeometric function with its parameters displaced by integers
is said to be contiguous to the original version, and the functions of
(\ref{eq:hgdef}) and~(\ref{eq:hypg}) are accordingly contiguous in a
generalized sense.  Any hypergeometric function with positive integral
differences between upper and lower parameters can be expressed as a
finite sum of hypergeometric functions of lower order, by what is now
called the Karlsson--Minton reduction formula~\cite{Karlsson71}.
Thus, the ${}_{r+\ell+1}F_{r+\ell}$ in~(\ref{eq:hypg}) can optionally
be written as a finite sum of ${}_{r+1}F_r$'s, though this fact will
not be exploited.

The key lemma used below is the following
(cf.~\cite[(5.7)]{Chen2005}).  Here, $\Delta(m;\mu)$ for $m\ge1$
abbreviates the $m$-parameter array
$\left(\frac{\mu}{m},\frac1m+\frac{\mu}{m},\dots,\frac{m-1}{m}+\frac{\mu}{m}\right)$.

\begin{lemma}
\label{lem:key}
  For\/ $l,m\ge1$, arbitrary parameter arrays\/ $(\alpha),(\beta)$ of
  lengths\/ $A,B$, and arbitrary\/~$a$ and\/ $x_0\neq0$, one has the identity
  \begin{equation}
  \begin{split}
    &{}_{l+m+A}F_B \left[
    \begin{array}{cc}
      \Delta(l+m;\,{a}), & (\alpha) \\
                         & (\beta)
    \end{array}
    \biggm|
    \frac{(l+m)^{l+m}}{l^l\;m^m}
    \,
    \frac{(-x/x_0)^l}{(1-x/x_0)^{l+m}}
    \right]
    \\[\jot]
    &\qquad=\left(1-{x}/{x_0}\right)^{{a}}
    F
    \left[
      \begin{array}{c:c|c}
        \begin{array}{c}
          \myatop{\vphantom{b}\textstyle a}{\vphantom{b}\textstyle\text{--}}
        \end{array}
        &
        \,
        R(n)
        &
        \,
        x/x_0
      \end{array}
      \right]
  \end{split}
  \end{equation}
  where
  \begin{displaymath}
    R(n) = {}_{l+m+A}F_B
    \left[
      \begin{array}{ccc}
        \Delta(l;\,-n), &  \Delta(m;\,n+{a}), & (\alpha)\\
        & & (\beta)
      \end{array}
      \biggm|
      1\,
      \right],
  \end{displaymath}
assuming the convergence of the series for the latter\/  ${}_{l+m+A}F_{B}(1)$.
\end{lemma}

\noindent
Only the case $l+m+A=B+1$ will be needed.  This is an identity of the
double-summation type: to prove it, one expands the hypergeometric
argument $\varphi(x)$ of the left-hand ${}_{l+m+A}F_B$ in a geometric
series, and converts the left side (multiplied by $(1-\nobreak
x/x_0)^{-{a}}$) to the right (multiplied by same) by interchanging the
order of the two summations.  It could be called classical; it was
stated by Bailey~\cite[\S\,4]{Bailey28}, and the $l=m=1$ case was
rediscovered by Chaundy and Rainville.  A~substantial generalization
was proved in~\cite{Chen2005}.  Special subcases of the $l=m=1$ case
are scattered in the literature; for details,
see~\cite[\S\,2.6]{Srivastava84}.

\section{Main theorems}
\label{sec:maintheorems}

The theorems can be thought of as being arranged in a $3\times3$
array.  Sections \ref{subsec:extensions}, \ref{subsec:gens},
and~\ref{subsec:extcompanions} contain the extensions of the classical
transformations of ${}_3F_2$ to itself, their generalizations
to~${}_4F_3$, and the extensions of the companion transformations of
${}_3F_2$ to~${}_4F_3$.  Each of these sections contains three
transformations: one quadratic and two cubic.

\subsection{Extended transformations of ${}_3F_2$}
\label{subsec:extensions}

The following theorems, indexed by $k\ge0$, reduce to Whipple's
quadratic transformation and Bailey's two cubic ones when $k=0$.  In
each, the left-hand ${}_3F_2$ has $S=\tfrac12+k$.

\begin{theorem}
\label{thm:A2}
  For all\/ $k\ge0$, one has the quadratic transformation
  \begin{align*}
    &{}_3F_2\left[
      \begin{array}{c}
        {\tfrac{a}2,\:\:\tfrac12+\tfrac{a}2,\:\:1-k+a-b-c}
        \\
          {1+a-b,\:\: 1+a-c}
      \end{array}
      \biggm|
      -\,\frac{4x}{(1-x)^2}
      \right]\\
    &\qquad=(1-x)^a\: {}_{3+2k}F_{2+2k}
    \left[
      \begin{array}{c:c|c}
        \begin{array}{ccc}
          a, & b, & c \\
          & 1+a-b, & 1+a-c
        \end{array}
        &
        \,
        Q_k^{(2)}(n)      
        &
        \,
        x
      \end{array}
      \right],
  \end{align*}
  where\/ $Q_k^{(2)}(n)=Q_k^{(2)}(n;a;b,c)$ is a degree-\/$2k$
  polynomial in\/~$n$ or a degree-\/$k$ one in\/
  $\lambda=\lambda(n;a)=n(n+a)$, the coordinate of a quadratic
  lattice, defined by
  \begin{displaymath}
  Q_k^{(2)}(n;\,a;\,b,c) = {}_3F_2
  \left[
    \begin{array}{c}
      {-n,\:n+a,\:-k} \\
      {b,\:c}
    \end{array}
  \right].
  \end{displaymath}
\end{theorem}

Here, the right-hand ${}_{3+2k}F_{2+2k}$ is well-poised for
all~$k\ge0$.  Owing to the $n\mapsto\allowbreak-n-\nobreak{a}$
invariance, the negated roots $\xi_1,\dots,\xi_{2k}$ of $Q_k^{(2)}$
are symmetric about $\xi=\tfrac{a}2$, and the lower parameters
$\xi_1,\dots,\xi_{2k}$ of the ${}_{3+2k}F_{2+2k}$ that are implicit in
this formula (recall~(\ref{eq:hypg})) can be permuted so that each
parameter-pair sums to~$1+\nobreak a$.

The $k=1$ case of this quadratic ${}_3F_2$ transformation, the first
to exhibit nonlinear parametric constraints, was discovered by
Niblett~\cite[(22)]{Niblett51}.  One finds
\begin{equation}
\label{eq:Q12}
Q_1^{(2)}(n;\,a;\,b,c) = 1 + \frac{\lambda}{bc} = \frac{n^2 + a n + bc}{bc},
\end{equation}
suggesting a subcase of interest: if $a=-b-c$, then
$Q_1^{(2)}(n)=(n-b)(n-c)/bc$ and the negated roots $\{\xi_1,\xi_2\}$
are $\{-b, -c\}$.  The resulting specialization is
\begin{equation}
\begin{split}
\label{eq:linconstraint}
  &{}_3F_2
  \left[
    \begin{array}{c}
      {-\tfrac{b}2-\tfrac{c}2,\:\:\tfrac12-\tfrac{b}2-\tfrac{c}2,\:\:-2b-2c}\\
      {1-2b-c,\:\:1-b-2c}
    \end{array}
    \biggm|
    -\,\frac{4x}{(1-x)^2}
    \right]
  \\
  &\qquad=(1-x)^{-b-c}\:
  {}_5F_4
  \left[
    \begin{array}{ccc:cc}
      -b-c, & b, & c, & 1-b, & 1-c \\
            &1-2b-c, & 1-b-2c, & -c, & -b
    \end{array}
    \biggm|
    x\,
    \right],
\end{split}
\end{equation}
in which the hypergeometric parameters are constrained linearly.  The
left-hand ${}_3F_2$ has $S=\frac12+k=\frac32$, and the right-hand
${}_5F_4$ is manifestly well-poised: the sum of each of its
parameter-pairs is $1-\nobreak b-\nobreak c$.  (Compare
\cite[(16)]{Niblett51}.)  Another notable $k=1$ subcase occurs when
$a^2=1+4bc$.  Then, the negated roots $\xi_1,\xi_2$ of $n^2+an+bc$
differ by unity, and the parameter-pairs
$\left[\begin{smallmatrix}1+\xi_1, & 1+\xi_2 \\ \xi_1, &
    \xi_2\end{smallmatrix}\right]$ can be merged into
$\left[\begin{smallmatrix}2+\xi_1\\ \xi_1\end{smallmatrix}\right]$,
reducing the right-hand ${}_5F_4$ to a~${}_4F_3$.

Other specializations of interest include the case
$c=\tfrac12+\tfrac{a}2$, when the transformation reduces to one of a
${}_2F_1$ with $S=\tfrac12+k$ to a well-poised ${}_{2+2k}F_{1+2k}$,
namely
\begin{equation}
  \begin{split}
    &{}_2F_1\left[
      \begin{array}{c}
        {\tfrac{a}2,\:\:\tfrac12-k+\tfrac{a}2-b}
        \\
          {1+a-b}
      \end{array}
      \biggm|
      -\,\frac{4x}{(1-x)^2}
      \right]
    \\
    &\qquad=(1-x)^a\: {}_{2+2k}F_{1+2k}
    \left[
      \begin{array}{c:c|c}
        \begin{array}{cc}
          a, & b \\
          & 1+a-b
        \end{array}
        &
        \,
        Q_k^{(2)}(n)      
        &
        \,
        x
      \end{array}
      \right],
  \end{split}
  \label{eq:lastmin2}
\end{equation}
where\/ $Q_k^{(2)}(n)\defeq Q_k^{(2)}(n;a;b,\tfrac12+\tfrac{a}2)$.
The $k=0$ subcase of~(\ref{eq:lastmin2}) is
classical~\cite[Thm.~3.1.1]{Andrews99}, but the $k>0$ subcases are
new.

By setting $x=-1$ in~(\ref{eq:lastmin2}), convergence of the series
being assumed, and evaluating the resulting ${}_2F_1(1)$ on the
left-hand side with the aid of Gauss's summation formula and the
duplication formula for the gamma function, one finds
\begin{equation}
\label{eq:reflastmin}
{}_{2+2k}F_{1+2k}
\left[
  \begin{array}{c:c|c}
    \begin{array}{cc}
      a, & b\\
      & 1+a-b
    \end{array}
    &
    \,
    Q_k^{(2)}(n)
    &
    \,
    -1
  \end{array}
\right]
=
\frac{(2k)!}{k!}
\,
\frac{\Gamma(1+a-b)\Gamma(1+k+\frac{a}2)}{\Gamma(1+2k+a)\Gamma(1+\frac{a}2-b)},
\end{equation}
where as before, $Q_k^{(2)}(n)\defeq
Q_k^{(2)}(n;a;b,\tfrac12+\tfrac{a}2)$.  Equation~(\ref{eq:reflastmin})
is an extension of Kummer's summation
formula~\cite[Cor.~3.1.2]{Andrews99} for a convergent, well-poised
${}_2F_1(-1)$, to which it reduces when $k=0$.  For all $k\ge0$, the
${}_{2+2k}F_{1+2k}(-1)$ series is well-poised and has $S=1-\nobreak
2k-\nobreak 2b$.  This is an extension of a type not previously
considered in the literature.

\begin{theorem}
\label{thm:A3}
For all\/ $k\ge0$, one has the first cubic transformation
\begin{displaymath}
  \begin{split}
    &{}_3F_2\left[
      \begin{array}{c}
        {\tfrac{a}3,\:\:\tfrac13+\tfrac{a}3,\:\:\tfrac23+\tfrac{a}3}
        \\
          {\tfrac34+\tfrac{k}2+\tfrac{a}2+\tfrac{b}2,\:\: \tfrac34+\tfrac{k}2+\tfrac{a}2-\tfrac{b}2}
      \end{array}
      \biggm|
      -\,\frac{27x}{(1-4x)^3}
      \right]\\
    &\quad=(1-4x)^a\: {}_{3+2k}F_{2+2k}
    \left[
      \begin{array}{c:c|c}
        \begin{array}{ccc}
          a, & \tfrac12-k-b, & \tfrac12-k+b \\
          & \tfrac34+\tfrac{k}2+\tfrac{a}2+\tfrac{b}2, & \tfrac34+\tfrac{k}2+\tfrac{a}2-\tfrac{b}2
        \end{array}
        &
        \,
        Q_k^{(3)}(n)      
        &
        \,
        x
        \end{array}
      \right],
  \end{split}
\end{displaymath}
where\/ $Q_k^{(3)}(n)=Q_k^{(3)}(n;a;b)$ is a degree-\/$2k$
polynomial in\/~$n$, equal to
\begin{equation}
  \frac{4^k(\tfrac14-\tfrac{k}2+\tfrac{b}2-\tfrac{n}2)_k   (\tfrac14-\tfrac{k}2-\tfrac{b}2-\tfrac{n}2)_k}{(\tfrac12+b)_k (\tfrac12-b)_k}\,
       {}_3F_2
       \left[
         {
           \myatop
               {-n,\:\frac{n}2+\frac{a}2,\:-k}
               {\tfrac14-\tfrac{k}2+\tfrac{b}2-\tfrac{n}2,\:\tfrac14-\tfrac{k}2-\tfrac{b}2-\tfrac{n}2}
         }
         \right].
       \nonumber
\end{equation}
\end{theorem}

Here, the right-hand ${}_{3+2k}F_{2+2k}$ is $(1,2)$-poised if $k=0$
(the classical case), but not otherwise.
This is illustrated by the $k=1$ case.  One finds
\begin{equation}
Q_1^{(3)}(n;\,a;\,b) = \frac{12n^2 + 4(1+2a)n + (1-4b^2)}{1-4b^2}
\end{equation}
(the denominator being required by the normalization
$Q_k^{(3)}(n=0)=1$; the subcase $b=\pm\tfrac12$ is singular).  From
this, the negated roots $\xi_1,\xi_2$ needed for the $k=1$ case can be
computed.  The resulting upper parameters $1+\nobreak\xi_1,\allowbreak
1+\nobreak\xi_2$ and lower ones $\xi_1,\xi_2$ implicit in the
right-hand ${}_{3+2k}F_{2+2k}={}_5F_4$ (recall~(\ref{eq:hypg})) do not
have the property that their sums (the lower ones being doubled) equal
$2+a$.

In the $k=1$ subcase with $a=-\frac12 \pm \frac{\sqrt3}2$,
$Q_1^{(3)}$~is proportional to $(n+\xi_1)(n+\xi_2)$ for
$\{\xi_1,\xi_2\}$ equal to
$\{\pm\frac{\sqrt{3}}{6}(1-2b),\pm\frac{\sqrt{3}}{6}(1+2b)\}$.  The
resulting specialization is
\begin{sizeequation}{\small}
\label{eq:linconstraint2}
\begin{split}
  &{}_3F_2
  \left[
    \begin{array}{c}
      {-\tfrac16\pm\tfrac{\sqrt{3}}6,\:\:\tfrac16\pm\tfrac{\sqrt{3}}6,\:\:\tfrac12\pm\tfrac{\sqrt{3}}6}\\
      {1\pm\tfrac{\sqrt{3}}4 + \tfrac{b}2,\:\:1\pm\tfrac{\sqrt{3}}4 - \tfrac{b}2}
    \end{array}
    \biggm|
    -\,\frac{27x}{(1-4x)^3}
    \right]
  =
  (1-4x)^{-\tfrac12\pm\tfrac{\sqrt{3}}2}\\
  &\:{}\times {}_5F_4
  \left[
    \begin{array}{ccc:cc}
      -\tfrac12\pm\tfrac{\sqrt{3}}2, & -\tfrac12-b, & -\tfrac12+b, & 1\pm\tfrac{\sqrt{3}}6(1-2b), & 1\pm\tfrac{\sqrt{3}}6(1+2b) \\
      & 1\pm\tfrac{\sqrt{3}}4 + \tfrac{b}2, & 1\pm\tfrac{\sqrt{3}}4 - \tfrac{b}2, & \pm\tfrac{\sqrt{3}}6(1+2b), & \pm\tfrac{\sqrt{3}}6(1-2b)
    \end{array}
    \biggm|
    x\,
    \right],
\end{split}
\end{sizeequation}
in which the parameters are constrained linearly.  It is analogous
to~(\ref{eq:linconstraint}).  The ${}_5F_4$ in this identity is
neither well-poised nor $(1,2)$-poised: to the right of the dashed
line, each parameter-pair sums to a constant (i.e.,
$1\pm\frac{\sqrt3}3$), but to the left, each sums to a constant (i.e.,
$\frac32\pm\frac{\sqrt3}2$) only if the lower member is doubled.

It is easily checked that when $\tfrac13\left(a+\tfrac12\right)^2 +
b^2 = 1$, the negated roots $\xi_1,\xi_2$ of $Q_1^{(3)}(n;\,a;\,b)$
differ by unity, i.e., $\xi_2=1+\xi_1$, allowing the right-hand
${}_{3+2k}F_{2+2k}={}_5F_4$ to be reduced to a~${}_4F_3$.  This
quadratic constraint on~$a,b$ (graphically, an ellipse) has the
parametrization $a=-\frac12+\sqrt{3}\,\sin\theta$, $b=\cos\theta$.  By
substituting into the $k=1$ case of the theorem, one
obtains~(\ref{eq:curious}).

Other specializations of interest include the case
$b=\tfrac16+k+\tfrac{a}3$, when the transformation formula reduces to
\begin{equation}
  \begin{split}
    &{}_2F_1\left[
      \begin{array}{c}
        {\tfrac{a}3,\:\:\tfrac13+\tfrac{a}3}
        \\[\jot]
          {\tfrac56+k+\tfrac{2a}3}
      \end{array}
      \biggm|
      -\,\frac{27x}{(1-4x)^3}
      \right]\\
    &\qquad=(1-4x)^a\: {}_{2+2k}F_{1+2k}
    \left[
      \begin{array}{c:c|c}
        \begin{array}{cc}
          a, & \tfrac13-2k-\tfrac{a}3 \\[\jot]
          & \tfrac56+k+\tfrac{2a}3
        \end{array}
        &
        \,
        Q_k^{(3)}(n)      
        &
        \,
        x
      \end{array}
      \right],
  \end{split}
  \label{eq:lastmin3}
\end{equation}
where\/ $Q_k^{(3)}(n)\defeq Q_k^{(3)}(n;a;b=\tfrac16+k+\tfrac{a}3)$.
The $k=0$ subcase of~(\ref{eq:lastmin3}) is a classically known cubic
transformation of a ${}_2F_1$ with $S=\frac12$ to a $(1,2)$-poised
${}_2F_1$, and is a specialization of Bailey's first cubic
transformation of~${}_3F_2$.  But the $k>0$ subcases are new.

\begin{theorem}
  \label{thm:A3p}
  For all\/ $k\ge0$, one has the second cubic transformation
  \begin{displaymath}
    \begin{split}
      &{}_3F_2\left[
        \begin{array}{c}
          {\tfrac{a}3,\:\:\tfrac13+\tfrac{a}3,\:\:\tfrac23+\tfrac{a}3}
          \\
            {\tfrac34+\tfrac{k}2+\tfrac{a}2+\tfrac{b}2,\:\: \tfrac34+\tfrac{k}2+\tfrac{a}2-\tfrac{b}2}
        \end{array}
        \biggm|
        \frac{27x^2}{(4-x)^3}
        \right]\\
      &\quad=\left(1-\tfrac{x}4\right)^a\, {}_{3+2k}F_{2+2k}
      \left[
        \begin{array}{c:c|c}
          \begin{array}{ccc}
            a, & \tfrac14 - \tfrac{k}2 + \tfrac{a}2 - \tfrac{b}2, & \tfrac14 - \tfrac{k}2 + \tfrac{a}2 + \tfrac{b}2 \\[\jot]
            & \tfrac12 + k + a + b, & \tfrac12 + k + a - b
          \end{array}
          &
          \,
          Q_k^{(3')}(n)      
          &
          \,
          x
        \end{array}
        \right],
    \end{split}
  \end{displaymath}
  where\/ $Q_k^{(3')}(n)=Q_k^{(3')}(n;a;b)$ is a degree-\/$2k$
  polynomial in\/~$n$, equal to
  \begin{displaymath}
    \begin{split}
      &\frac{(\tfrac34-\tfrac{k}2-\tfrac{a}2+\tfrac{b}2-n)_k(\tfrac34-\tfrac{k}2-\tfrac{a}2-\tfrac{b}2-n)_k}{(\tfrac34-\tfrac{k}2-\tfrac{a}2+\tfrac{b}2)_k(\tfrac34-\tfrac{k}2-\tfrac{a}2-\tfrac{b}2)_k}
      \\
      &
      \qquad{}\times
            {}_3F_2
            \left[
              {
                \myatop
                    {-\frac{n}2,\:-\frac{n}2+\frac{1}2,\:-k}
                    {\tfrac34-\tfrac{k}2-\tfrac{a}2+\tfrac{b}2-{n},\:\tfrac34-\tfrac{k}2-\tfrac{a}2-\tfrac{b}2-{n}}
              }
              \right].
            \nonumber
    \end{split}
  \end{displaymath}
\end{theorem}

Here, the right-hand ${}_{3+2k}F_{2+2k}$ is $(2,1)$-poised if $k=0$
(the classical case), but not otherwise.  The polynomials $Q_k^{(3')}$
differ from the~$Q_k^{(3)}$; for instance,
\begin{equation}
  Q_1^{(3')}(n;\,a;\,b) =
  \frac{ 12n^2 - 4(1-4a)n + (1-2a-2b)(1-2a+2b) }{(1-2a-2b)(1-2a+2b)}.
\end{equation}
As with Theorem~\ref{thm:A3}, there are interesting specializations.

\subsection{Generalizations to ${}_4F_3$}
\label{subsec:gens}

Each left-hand ${}_4F_3$ in the following theorems has $S=\tfrac12$
and contains a parameter-pair $\left[\begin{smallmatrix}k+d
  \\ d\end{smallmatrix}\right]$, where $d$~is an additional free
  parameter.  These identities reduce to Whipple's and Bailey's
  classical transformations when $k=0$, and to the extensions
  of~\S\,\ref{subsec:extensions} when~$d\to\infty$.  It should be
  noted that by the Karlsson--Minton reduction
  formula~\cite{Karlsson71}, any ${}_4F_3$ with a
  parameter-pair~$\left[\begin{smallmatrix}k+d
    \\ d\end{smallmatrix}\right]$ can be written as a sum of
    $1+\nobreak k$ functions of the ${}_3F_2$ type.

\begin{theorem}
\label{thm:B2}
  For all\/ $k\ge0$, one has the quadratic transformation
  \begin{equation}
  \label{eq:lastminute}
  \begin{split}
      &{}_4F_3\left[
        \begin{array}{c:c}
          {\tfrac{a}2,\:\:\tfrac12+\tfrac{a}2,\:\:1-k+a-b-c,} & {k+d} \\
          {1+a-b,\:\: 1+a-c,} & {d}
        \end{array}
        \biggm|
        -\,\frac{4x}{(1-x)^2}
        \right]\\
      &\qquad=(1-x)^a\: {}_{3+2k}F_{2+2k}
      \left[
        \begin{array}{c:c|c}
          \begin{array}{ccc}
            a, & b, & c \\
            & 1+a-b, & 1+a-c
          \end{array}
          &
          \,
          Q_k^{(2)}(n)      
          &
          \,
          x
          \end{array}
        \right],
  \end{split}
  \end{equation}
  where\/ $Q_k^{(2)}(n)=Q_k^{(2)}(n;a;b,c,d)$ is a degree-\/$2k$
  polynomial in\/~$n$ or a degree-\/$k$ one in\/
  $\lambda=\lambda(n;a)=n(n+a)$, the coordinate of a quadratic
  lattice, defined by
  \begin{displaymath}
    Q_k^{(2)}(n;\,a;\,b,c,d) = {}_4F_3
    \left[
      \begin{array}{c:c}
        {-n,\:n+a,\:-k,} & {k-1-a+b+c+d} \\
        {b,\:c,} & {d}
      \end{array}
      \right].
  \end{displaymath}
\end{theorem}

Here, the right-hand ${}_{3+2k}F_{2+2k}$ is well-poised for all
$k\ge0$, as in Theorem~\ref{thm:A2}.  The four-parameter
$Q_k^{(2)}(n)$ is essentially a Racah
polynomial~\cite[\S\,9.2]{Koekoek2010}, just as the three-parameter
one in Theorem~\ref{thm:A2} was a dual Hahn polynomial.  For instance,
\begin{align}
\label{eq:Q1}
Q_1^{(2)}(n;\,a;\,b,c,d) &= 1 + \frac{(b+c+d-a)\lambda}{bcd}\\
 &= \frac{(b+c+d-a)n^2 + a(b+c+d-a)\,n + bcd}{bcd},\nonumber
\end{align}
the $d\to\infty$ limit of which is the $Q_1^{(2)}(n;a,b,c)$
of~(\ref{eq:Q12}).  Owing to the $n\mapsto-n-a$ invariance, the
negated roots $\xi_1,\dots,\xi_{2k}$ of~$Q^{(2)}_k$ are symmetric
about $\xi=\nobreak\tfrac{a}2$.

Specializations of interest include the choice $c=\tfrac12+\tfrac{a}2$,
which leads to
\begin{equation}
\begin{split}
\label{eq:rrplus}
  &{}_3F_2\left[
    \begin{array}{c:c}
      {\tfrac{a}2,\:\:\tfrac12-k+\tfrac{a}2-b,} & {k+d} \\
      {1+a-b,} & {d}
    \end{array}
    \biggm|
    -\,\frac{4x}{(1-x)^2}
    \right]\\
  &\qquad=(1-x)^a\: {}_{2+2k}F_{1+2k}
  \left[
    \begin{array}{c:c|c}
      \begin{array}{ccc}
        a, & b \\
        & 1+a-b
      \end{array}
      &
      \,
      Q_k^{(2)}(n)      
      &
      \,
      x
      \end{array}
    \right],
\end{split}
\end{equation}
where $Q_k^{(2)}(n)\defeq Q_k^{(2)}(n;a;b,\tfrac12+\tfrac{a}2;d)$.
The $k=1$ cases of (\ref{eq:lastminute}) and~(\ref{eq:rrplus}) are
known (see \cite[Thm.~1]{Wang2013}, resp.\ \cite[(3.1)]{Rakha2011}).
It must be mentioned that other transformations of a ${}_3F_2$ with a
parameter-pair~$\left[\begin{smallmatrix}1+d
  \\ d\end{smallmatrix}\right]$ to a~${}_4F_3$ have been found (see
  \cite[\S\,6]{Miller2013} and~\cite{Rakha2009}).  The others
  have lifting functions $\varphi(x)$ equal to $\tfrac{x^2}{(2-x)^2}$, 
  $\tfrac{4x}{(1+x)^2}$, and $4x(1-\nobreak x)$.

\begin{theorem}
\label{thm:B3}
For all\/ $k\ge0$, one has the first cubic transformation
\begin{displaymath}
  \begin{split}
    &{}_4F_3\left[
      \begin{array}{c:c}
        {\tfrac{a}3,\:\:\tfrac13+\tfrac{a}3,\:\:\tfrac23+\tfrac{a}3,} & k+d \\
        {\tfrac34+\tfrac{k}2+\tfrac{a}2+\tfrac{b}2,\:\: \tfrac34+\tfrac{k}2+\tfrac{a}2-\tfrac{b}2,} & d
      \end{array}
      \biggm|
      -\,\frac{27x}{(1-4x)^3}
      \right]\\
    &\quad=(1-4x)^a\: {}_{3+3k}F_{2+3k}
    \left[
      \begin{array}{c:c|c}
        \begin{array}{ccc}
          a, & \tfrac12-k-b, & \tfrac12-k+b \\
          & \tfrac34+\tfrac{k}2+\tfrac{a}2+\tfrac{b}2, & \tfrac34+\tfrac{k}2+\tfrac{a}2-\tfrac{b}2
        \end{array}
        &
        \,
        Q_k^{(3)}(n)      
        &
        \,
        x
      \end{array}
      \right],
  \end{split}
\end{displaymath}
where\/ $Q_k^{(3)}(n)=Q_k^{(3)}(n;a;b;d)$ is a degree-\/$3k$
polynomial in\/~$n$, defined as in Theorem\/~{\rm\ref{thm:A3}} but with
the\/ ${}_3F_2(1)$ in the definition extended to
\begin{displaymath}
  {}_4F_3
  \left[
    \begin{array}{c:c}
      {-n,\:\frac{n}2+\frac{a}2,\:-k,} & -\tfrac{n}2-\tfrac{a}2-\tfrac{1}2+d\\
      {\tfrac14-\tfrac{k}2+\tfrac{b}2-\tfrac{n}2,\:\tfrac14-\tfrac{k}2-\tfrac{b}2-\tfrac{n}2,} & d
    \end{array}
    \right].
\end{displaymath}
\end{theorem}

\begin{theorem}
  \label{thm:B3p}
  For all\/ $k\ge0$, one has the second cubic transformation
  \begin{displaymath}
    \begin{split}
      &{}_4F_3\left[
        \begin{array}{c:c}
          {\tfrac{a}3,\:\:\tfrac13+\tfrac{a}3,\:\:\tfrac23+\tfrac{a}3,} & k+d
          \\
            {\tfrac34+\tfrac{k}2+\tfrac{a}2+\tfrac{b}2,\:\: \tfrac34+\tfrac{k}2+\tfrac{a}2-\tfrac{b}2,} & d
        \end{array}
        \biggm|
        \frac{27x^2}{(4-x)^3}
        \right]\\
      &\quad=\left(1-\tfrac{x}4\right)^a\, {}_{3+3k}F_{2+3k}
      \left[
        \begin{array}{c:c|c}
          \begin{array}{ccc}
            a, & \tfrac14 - \tfrac{k}2 + \tfrac{a}2 - \tfrac{b}2, & \tfrac14 - \tfrac{k}2 + \tfrac{a}2 + \tfrac{b}2 \\
            & \tfrac12 + k + a + b, & \tfrac12 + k + a - b
          \end{array}
          &
          \,
          Q_k^{(3')}(n)      
          &
          \,
          x
        \end{array}
        \right],
    \end{split}
  \end{displaymath}
  where\/ $Q_k^{(3')}(n)=Q_k^{(3')}(n;a;b;d)$ is a degree-\/$3k$
  polynomial in\/~$n$, defined as in Theorem\/~{\rm\ref{thm:A3p}} but with
  the\/ ${}_3F_2(1)$ in the definition extended to
  \begin{displaymath}
    {}_4F_3
    \left[
      \begin{array}{c:c}
        {-\frac{n}2,\:-\frac{n}2+\frac{1}2,\:-k,} & -n-a+d \\
        {\tfrac34-\tfrac{k}2-\tfrac{a}2+\tfrac{b}2-{n},\:\tfrac34-\tfrac{k}2-\tfrac{a}2-\tfrac{b}2-{n},} & d
      \end{array}
      \right].
  \end{displaymath}
\end{theorem}

\subsection{Extended companion transformations of ${}_3F_2$}
\label{subsec:extcompanions}

The following theorems, indexed by $k\ge0$, reduce to the companions
of Whipple's quadratic transformation and Bailey's two cubic ones when
$k=0$.  (For the companions, see \cite[p.~97,\ Example~6]{Bailey35}
and \cite[(4.1),(5.4),(5.7)]{Gessel82}.)  In each, the left-hand
${}_3F_2$ has $S=-\tfrac12-k$.

\begin{theorem}
  \label{thm:C2}
  For all\/ $k\ge0$, one has the quadratic transformation
  \begin{displaymath}
    \begin{split}
      &{}_3F_2\left[
        \begin{array}{c}
          {\tfrac12 + k + \tfrac{a}2,\:\:1+k+\tfrac{a}2,\:\:1-k+a-b-c}
          \\
            {1+a-b,\:\: 1+a-c}
        \end{array}
        \biggm|
        -\,\frac{4x}{(1-x)^2}
        \right]= (1+x)^{-1-2k} \\
      &\qquad{}\times(1-x)^{1+2k+a}\: {}_{4+4k}F_{3+4k}
      \left[
        \begin{array}{c:c|c}
          \begin{array}{ccc}
            a, & b, & c \\
            & 1+a-b, & 1+a-c
          \end{array}
          &
          \,
            {\mathbf{Q}}_k^{(2)}(n)      
            &
            \,
            x
        \end{array}
        \right],
    \end{split}
  \end{displaymath}
where ${\mathbf{Q}}^{(2)}_k(n) = {\mathbf{Q}}_k^{(2)}(n;a;b,c)$ is a
degree-\/$(1+\nobreak 4k)$ polynomial in\/ $n$, equal to\/
$1+\nobreak\tfrac{2n}a$ times\/ $\hat{\mathbf{Q}}^{(2)}_k(n) =
\hat{\mathbf{Q}}^{(2)}_k(n;a;b,c)$, which is a degree-\/$4k$
polynomial in\/ $n$ or a degree-\/$2k$ one in\/
$\lambda=\lambda(n;a)=n(n+a)$, the coordinate of a quadratic lattice,
determined by\/ $\hat{\mathbf{Q}}^{(2)}_0\equiv1$ and the\/
$k$\nobreakdash-raising relation
\begin{align*}
& (k+\tfrac{a}2)(1+a)bc\, (n+\tfrac{a}2)\, \hat{\mathbf{Q}}_k^{(2)}(n) \\
&\qquad=(n+k+\tfrac{a}2)(n+a)(n+b)(n+c)(n+\tfrac{1}2+\tfrac{a}2)\, \hat{\mathbf{Q}}_{k-1,+}^{(2)}(n) \nonumber\\
&\qquad\quad{}+ (n-k+\tfrac{a}2)n(n+a-b)(n+a-c)(n-\tfrac{1}2+\tfrac{a}2)\, \hat{\mathbf{Q}}_{k-1,+}^{(2)}(n-1) ,  \nonumber
\end{align*}
with\/ $\hat{\mathbf{Q}}_{k-1,+}^{(2)}(n) \defeq
\hat{\mathbf{Q}}^{(2)}_{k-1}(n;a+1,b+1,c+1)$.
\end{theorem}

Here, the right-hand ${}_{4+4k}F_{3+3k}$ is very well poised for all
$k\ge0$, because one negated root is $\xi_1=\tfrac{a}2$, coming from
the factor $1+\frac{2n}a$, and the remaining ones
$\xi_2,\dots,\xi_{1+4k}$ are symmetric about $\xi=\tfrac{a}2$, as
the recurrence for $\hat{\mathbf{Q}}_k^{(2)}$ is invariant under
$n\mapsto-n-a$. An example of $\hat{\mathbf{Q}}_k^{(2)}$ being of
degree $2k$ in $\lambda=n(n+a)$ is
\begin{equation}
  \hat{\mathbf{Q}}_1^{(2)}(n;\,a;\,b,c)=1+
  \frac{\lambda\left[4\lambda+(a-1)(a-2)+(2b+3)(2c+3)-9\right]}{(a+1)(a+2)bc}.
\end{equation}

Specializations of interest include the case $c=\tfrac12+\tfrac{a}2$,
when the right-hand ${}_{4+4k}F_{3+4k}$ reduces to a ${}_{3+4k}
F_{2+4k}$, and $c=\tfrac12-k+\tfrac{a}2$ and $c=-k+\tfrac{a}2$, when
the left-hand ${}_3F_2$ reduces to a~${}_2F_1$.  One can show from the
raising relation that, e.g.,
\begin{equation}
  \hat{\mathbf{Q}}^{(2)}_k(n;\,a;\,b,\tfrac12-k+\tfrac{a}2) =
  \frac{(\tfrac12+k+\tfrac{a}2)_n}{(\tfrac12-k+\tfrac{a}2)_n}
  \:\,{}_3F_2\left[
    \begin{array}{c}
      -n,\:\: n+a, \:\: -k \\
      b, \:\: 1+\tfrac{a}2
    \end{array}
    \right],
\end{equation}
$2k$ of the $4k$~negated roots of which are
$\tfrac12-k+\tfrac{a}2,\dots, -\tfrac12+k+\tfrac{a}2$.  But for
general parameter choices, a hypergeometric representation of
$\hat{\mathbf{Q}}^{(2)}_k(n;a;b,c)$ is lacking.

\begin{theorem}
  \label{thm:C3}
  For all\/ $k\ge0$, one has the first cubic transformation
  {
    \arraycolsep0.25em
    \begin{sizedisplaymath}{\small}
      \begin{split}
        &{}_3F_2\left[
          \begin{array}{c}
            {\tfrac13 + \tfrac{2k}3 + \tfrac{a}3,\:\:\tfrac23+ \tfrac{2k}3+\tfrac{a}3, \:\:1 + \tfrac{2k}3 + \tfrac{a}3}
            \\
              {\tfrac34+\tfrac{k}2+\tfrac{a}2+\tfrac{b}2,\:\: \tfrac34+\tfrac{k}2+\tfrac{a}2-\tfrac{b}2}
          \end{array}
          \biggm|
          -\,\frac{27x}{(1-4x)^3}
          \right]= 
        (1+8x)^{-1-2k}
        \\
        &\qquad{}\times(1-4x)^{1+2k+a}\, {}_{4+4k}F_{3+4k}\left[
          \begin{array}{c:c|c}
            \begin{array}{ccc}
              a, & \tfrac12-k-b, & \tfrac12-k+b \\
              & \tfrac34+\tfrac{k}2+\tfrac{a}2+\tfrac{b}2, & \tfrac34+\tfrac{k}2+\tfrac{a}2-\tfrac{b}2
            \end{array}
            &
            \,
              {\mathbf{Q}}_k^{(3)}(n)
              &
              \,
              x
          \end{array}
          \right],
      \end{split}
  \end{sizedisplaymath}}where
  ${\mathbf{Q}}_k^{(3)}(n) = {\mathbf{Q}}_k^{(3)}(n;a;b)$ is a
  degree-\/$(1+\nobreak 4k)$ polynomial in\/ $n$, determined by\/
  ${\mathbf{Q}}^{(3)}_0 = 1+\nobreak\tfrac{3n}a$ and the\/ $k$\nobreakdash-raising relation
  \begin{align*}
    & a(\tfrac12-k-b)(\tfrac12-k+b)\, {\mathbf{Q}}_k^{(3)}(n) \\
    &\quad=\left(\tfrac{3n+2k+a}{2k+a} \right)(n+a)(n+\tfrac12-k-b)(n+\tfrac12-k+b)\, {\mathbf{Q}}_{k-1,+}^{(3)}(n) \nonumber\\
    &\quad\quad{}+ 8\left(\tfrac{3n-4k+a}{2k+a} \right)   n(n-\tfrac14+\tfrac{k}2+\tfrac{a}2+\tfrac{b}2)(n-\tfrac14+\tfrac{k}2+\tfrac{a}2-\tfrac{b}2)\, {\mathbf{Q}}_{k-1,+}^{(3)}(n-1) ,  \nonumber
  \end{align*}
  with\/ ${\mathbf{Q}}_{k-1,+}^{(3)}(n) \defeq
  {\mathbf{Q}}^{(3)}_{k-1}(n;a+1,b)$.
\end{theorem}

\begin{theorem}
  \label{thm:C3p}
  For all\/ $k\ge0$, one has the second cubic transformation
  {
    \arraycolsep0.2em
    \begin{sizedisplaymath}{\small}
      \begin{split}
        &{}_3F_2\left[
          \begin{array}{c}
            {\tfrac13 + \tfrac{2k}3 + \tfrac{a}3,\:\:\tfrac23+ \tfrac{2k}3+\tfrac{a}3, \:\:1 + \tfrac{2k}3 + \tfrac{a}3}
            \\
              {\tfrac34+\tfrac{k}2+\tfrac{a}2+\tfrac{b}2,\:\: \tfrac34+\tfrac{k}2+\tfrac{a}2-\tfrac{b}2}
          \end{array}
          \biggm|
          \frac{27x^2}{(4-x)^3}
          \right]= 
        \left(1+\tfrac{x}8\right)^{-1-2k}
        \\
        &\qquad{}\times\left(1-\tfrac{x}4\right)^{1+2k+a}\, {}_{4+4k}F_{3+4k}\left[
          \begin{array}{c:c|c}
            \begin{array}{ccc}
              a, & \tfrac14-\tfrac{k}2+\tfrac{a}2-\tfrac{b}2,  &  \tfrac14-\tfrac{k}2+\tfrac{a}2+\tfrac{b}2 \\
              & \tfrac12+k+a+b, & \tfrac12+k+a-b
            \end{array}
            &
            \,
                {\mathbf{Q}}_k^{(3')}(n)
                &
                \,
                x
          \end{array}
          \right], \nonumber
      \end{split}
  \end{sizedisplaymath}}where 
  ${\mathbf{Q}}_k^{(3')}(n) = {\mathbf{Q}}_k^{(3')}(n;a;b)$ is a
  degree-\/$(1+\nobreak 4k)$ polynomial in\/ $n$, determined by\/
  ${\mathbf{Q}}^{(3')}_0 = 1+\nobreak\tfrac{3n}{2a}$ and the\/ $k$\nobreakdash-raising relation
  \begin{align*}
    & 2a(\tfrac14-\tfrac{k}2+\tfrac{a}2-\tfrac{b}2)(\tfrac14-\tfrac{k}2+\tfrac{a}2+\tfrac{b}2)\, {\mathbf{Q}}_k^{(3')}(n)  \\
    &\quad=\left(\tfrac{3n+4k+2a}{2k+a} \right)(n+a)(n+\tfrac14-\tfrac{k}2+\tfrac{a}2-\tfrac{b}2)(n+\tfrac14-\tfrac{k}2+\tfrac{a}2+\tfrac{b}2)\, {\mathbf{Q}}_{k-1,+}^{(3')}(n) \nonumber\\
    &\quad\quad{}+ \frac18\left(\tfrac{3n-2k+2a}{2k+a} \right)   n(n-\tfrac12+k+a+b)(n-\tfrac12+k+a-b)\, {\mathbf{Q}}_{k-1,+}^{(3')}(n-1) ,  \nonumber
  \end{align*}
  with\/ ${\mathbf{Q}}_{k-1,+}^{(3')}(n) \defeq
  {\mathbf{Q}}^{(3')}_{k-1}(n;a+1,b)$.
\end{theorem}

\section{Proofs}
\label{sec:proofs}

The following are the proofs of the first six theorems
of~\S\,\ref{sec:maintheorems}, those of the final three being deferred
to the next section.  The proofs of the first three employ the
Sheppard--Andersen transformation of terminating ${}_3F_2(1)$'s, which
is \cite[Cor.~3.3.4]{Andrews99}
\begin{align}
\label{eq:sheppard}
  {}_3F_2\left[
  \begin{array}{ccc}
    -n, & A, & B \\ & D, & E
  \end{array}
  \right]
  &=
  \left[
    \begin{array}{c}
      D-A,\:\:E-A\\
      D,\:\:E
    \end{array}
    \right]_n
  \\
  &\qquad{}\times  {}_3F_2
    \left[
    \begin{array}{c}
      -n, \:\: A, \:\: 1-S \\
      1+A-D-n, \:\: 1+A-E-n
    \end{array}
    \right],\nonumber
\end{align}
where $S=n-A-B+D+E$ is the parametric excess of the left-hand
${}_3F_2(1)$.  The notation $\left[\begin{smallmatrix} (\alpha)
    \\ (\beta) \end{smallmatrix}\right]_n$ signifies $
\tfrac{((\alpha))_n}{((\beta))_n}= \tfrac{\prod_i
  (\alpha_i)_n}{\prod_i (\beta_i)_n} $.

The formula~(\ref{eq:sheppard}) specializes when $S=1$ to the
Pfaff--Saalsch\"utz formula for the sum of any
$1$\nobreakdash-balanced terminating ${}_3F_2(1)$ series.  Also,
(\ref{eq:sheppard}) extends to Whipple's transformation of
$1$\nobreakdash-balanced terminating ${}_4F_3(1)$'s, which
is~\cite[Thm.~3.3.3]{Andrews99}
\begin{align}
\label{eq:whipple}
  {}_4F_3\left[
  \begin{array}{ccc:c}
    -n, & A, & B, & C \\ & D, & E, & F
  \end{array}
  \right]
  &=
  \left[
    \begin{array}{c}
      D-A,\:\:E-A\\
      D,\:\:E
    \end{array}
    \right]_n
  \\
  &\qquad{}\times
    {}_4F_3
    \left[
    \begin{array}{c:c}
      -n, \:\: A, \:\: F-B, & F-C\\
      1+A-D-n, \:\: 1+A-E-n, & F
    \end{array}
    \right].\nonumber
\end{align}
It is assumed in~(\ref{eq:whipple}) that the parametric excess of the
left-hand ${}_4F_3(1)$, which is $n-\nobreak A-B-C+D+E+\nobreak F$,
equals unity.  Equation~(\ref{eq:whipple}) can be deduced from Euler's
transformation of~${}_2F_1$, and (\ref{eq:sheppard}) comes
from~(\ref{eq:whipple}) by taking $C,F\to\infty$ with
$F-\nobreak{C}=\textrm{const}$.

The quadratic identity of Theorem~\ref{thm:A2} and the cubic ones of
Theorems \ref{thm:A3} and \ref{thm:A3p} follow respectively from the
$(l,m)=\allowbreak(1,1)$, $(1,2)$, and~$(2,1)$ cases of
Lemma~\ref{lem:key}, provided that the parameter arrays
$(\alpha),(\beta)$ are defined by
\begin{align*}
&(l,m)=(1,1):& \qquad  (\alpha)=(1-k+a-b-c), \quad (\beta)=(1+a-b,\, 1+a-c), \\
&(l,m)=(1,2): &\qquad  (\alpha)=(\textrm{--}), \quad
(\beta) = (\tfrac34+\tfrac{k}2+\tfrac{a}2+\tfrac{b}2,\,
\tfrac34+\tfrac{k}2+\tfrac{a}2-\tfrac{b}2),\\
&(l,m)=(2,1): &\qquad  (\alpha)=(\textrm{--}), \quad
(\beta) = (\tfrac34+\tfrac{k}2+\tfrac{a}2+\tfrac{b}2,\,
\tfrac34+\tfrac{k}2+\tfrac{a}2-\tfrac{b}2),
\end{align*}
with $x_0$ taken respectively to equal $1,\tfrac14,4$.  The
$(l,m)=\allowbreak (1,1)$, $(1,2)$ cases of the lemma can then be
written as
\begin{subequations}
  \begin{align}
    &
    \begin{aligned}
      &
      {}_3F_2
      \left[
        \begin{array}{c}
          \tfrac{a}2, \:\: \tfrac12 + \tfrac{a}2, \:\: 1-k+a-b-c \\
          1+a-b, \:\: 1+a-c
        \end{array}
        \biggm|
        -\, \frac{4x}{(1-x)^2}
        \right]
      \\
      &\qquad\qquad=(1-x)^a\:{}_3F_2
      \left[
        \begin{array}{c:c|c}
          \begin{array}{c}
            a\\
            \textrm{--}
          \end{array}
          &
          \,
          R^{(2)}(n)
          &
          \,
          x
          \end{array}
        \right],
    \end{aligned}
    \label{eq:D2}
    \\[\jot]
    &
    \begin{aligned}
      &{}_3F_2
      \left[
        \begin{array}{c}
          \tfrac{a}3, \:\: \tfrac13 + \tfrac{a}3, \:\: \tfrac{2}3 + \tfrac{a}3\\
          \tfrac34+\tfrac{k}2+\tfrac{a}2+\tfrac{b}2, \:\:
          \tfrac34+\tfrac{k}2+\tfrac{a}2-\tfrac{b}2
        \end{array}
        \biggm|
        -\, \frac{27x}{(1-4x)^3}
        \right]
      \\
      &\qquad\qquad=(1-4x)^{a}\:{}_3F_2
      \left[
        \begin{array}{c:c|c}
          \begin{array}{c}
            a \\ \textrm{--}
          \end{array}
          &
          \,
          R^{(3)}(n)
          &
          \,
          4x
        \end{array}
        \right],
    \end{aligned}
    \label{eq:D3}
  \end{align}
\end{subequations}
\noindent
where each ${}_3F_2$ has $S=1+k$ and each of $R^{(2)},R^{(3)}$ is a
${}_3F_2(1)$ with $S=1+\nobreak k$, i.e.,
\begin{subequations}
\begin{sizealign}{\normalsize}
\label{eq:E2}
  R^{(2)}(n) &= {}_3F_2
    \left[
      \begin{array}{c}
        -n,\:\: n+a, \:\: 1-k+a-b-c\\
        1+a-b, \:\: 1+a-c
      \end{array} 
      \right]
    \\
    &=
    \left[
      \begin{array}{cc}
        b, & c\\
        1+a-b, & 1+a-c
      \end{array}
      \right]_n
    \: 
      {}_3F_2
      \left[
        \begin{array}{c}
          -n, \:\: n+a, \:\: -k\\
          b, \:\: c
        \end{array}
        \right],
      \nonumber\\[\jot]
      \label{eq:E3}
      R^{(3)}(n) &= {}_3F_2
      \left[
        \begin{array}{c}
          -n, \:\: \tfrac{n}2 +\tfrac{a}2, \:\: \tfrac{n}2 + \tfrac{1}2 + \tfrac{a}2\\
        \tfrac34+\tfrac{k}2+\tfrac{a}2+\tfrac{b}2, \:\:
        \tfrac34+\tfrac{k}2+\tfrac{a}2-\tfrac{b}2
        \end{array}
        \right]
      \\
      &=
      \left[
        \begin{array}{cc}
          \tfrac34 + \tfrac{k}2 + \tfrac{b}2 - \tfrac{n}2, & 
          \tfrac34 + \tfrac{k}2 - \tfrac{b}2 - \tfrac{n}2 \\
          \tfrac34 + \tfrac{k}2 + \tfrac{a}2 + \tfrac{b}2, &
          \tfrac34 + \tfrac{k}2 + \tfrac{a}2 - \tfrac{b}2
        \end{array}
        \right]_n
      \nonumber
      \\
      &\qquad{}\times
      {}_3F_2\left[
        \begin{array}{c}
          -n, \:\: \tfrac{n}2 + \tfrac{a}2, \:\: -k \\
          \tfrac14 - \tfrac{k}2 + \tfrac{b}2 - \tfrac{n}2, \:\:
          \tfrac14 - \tfrac{k}2 - \tfrac{b}2 - \tfrac{n}2
        \end{array}
        \right]. \nonumber
\end{sizealign}
\end{subequations}
The second expressions for $R^{(2)}(n)$, $R^{(3)}(n)$ are obtained by
applying the transformation~(\ref{eq:sheppard}).  The prefactor in the
second expression in~(\ref{eq:E3}) equals
\begin{sizeequation}{\small}
  \label{eq:redone}
  4^{-n}\left[
    \begin{array}{cc}
    \tfrac12-k-b, & \tfrac12-k+b \\
    \tfrac34+\tfrac{k}2+\tfrac{a}2+\tfrac{b}2, &
    \tfrac34+\tfrac{k}2+\tfrac{a}2-\tfrac{b}2
    \end{array}
    \right]_n
  \left[
  \frac{4^k(\tfrac14-\tfrac{k}2+\tfrac{b}2-\tfrac{n}2)_k   (\tfrac14-\tfrac{k}2-\tfrac{b}2-\tfrac{n}2)_k}{(\tfrac12+b)_k (\tfrac12-b)_k}
  \right]
\end{sizeequation}
by elementary Pochhammer-symbol manipulations.

Substituting (\ref{eq:E2}),(\ref{eq:E3}) [with~(\ref{eq:redone})] into
(\ref{eq:D2}),(\ref{eq:D3}) immediately yields the identities of
Theorems \ref{thm:A2} and~\ref{thm:A3}.  The derivation of the cubic
transformation in Theorem~\ref{thm:A3p} from the $(l,m)=(2,1)$ case of
the lemma proceeds similarly, with a minor difference: its
even\nobreakdash-$n$ and odd\nobreakdash-$n$ subcases must be treated
separately.

The proofs of the ${}_4F_3$ transformations in Theorems \ref{thm:B2},
\ref{thm:B3}, and~\ref{thm:B3p} are identical to the preceding three,
except that $(\alpha),(\beta)$ include $k+\nobreak{d},d$ as respective
additional parameters, each of $R^{(2)},R^{(3)}$ is a ${}_4F_3(1)$
with $S=1$, and Whipple's transformation~(\ref{eq:whipple}) is used
instead of the Sheppard--Andersen transformation~(\ref{eq:sheppard}).

It is worth recalling that there are no close analogues of the
Sheppard--Andersen and Whipple transformations for terminating
hypergeometric series of higher order than ${}_3F_2(1)$
and~${}_4F_3(1)$.  The known transformations of terminating
${}_7F_6(1)$ series exist because of Whipple's formula relating
certain ${}_7F_6(1)$'s and ${}_4F_3(1)$'s.  And although
transformations of terminating ${}_9F_8(1)$ series are known, the
series must satisfy restrictive conditions (e.g., they must be very
well poised as~well as $2$\nobreakdash-balanced).

\section{The polynomials $Q_k$ and $\mathbf{Q}_k$: Raising relations}
\label{sec:polys}

Each of the polynomials $Q_k^{(2)},Q_k^{(3)},Q_k^{(3')}$ and
${\mathbf{Q}}_k^{(2)},{\mathbf{Q}}_k^{(3)},{\mathbf{Q}}_k^{(3')}$ in
the transformations of~\S\,\ref{sec:maintheorems} satisfies a
recurrence on~$k$, to be called a $k$\nobreakdash-raising relation.
The~$Q_k$, appearing in the six transformations of
\S\S\,\ref{subsec:extensions} and~\ref{subsec:gens}, have
hypergeometric representations from which recurrences can be deduced.
But in all cases it is easier to go directly from a hypothesized
transformation $F(t) = \mathcal{A}(x)\tilde F(x)$, based on a lifting
function $t=\varphi(x)$, to the corresponding recurrence.  It will now
be shown how this can be done.  The $k$\nobreakdash-raising relations
satisfied by the polynomials $Q_k$ are of some importance, but those
satisfied by the polynomials ${\mathbf{Q}}_k$
of~\S\,\ref{subsec:extcompanions}, deduced from the statements of
Theorems \ref{thm:C2}, \ref{thm:C3}, and~\ref{thm:C3p} in that
subsection, are of particular importance: they constitute a proof of
those theorems, by induction on~$k$.  This is because the classical
($k=0$) cases of those theorems have long been known.

Suppose $t=\varphi(x)$ and that $F(t),\tilde F(x)$ are hypergeometric
functions.  Define $\vartheta=t\tfrac{{\rm d}}{{\rm d}t}$ and
$\delta=x\tfrac{{\rm d}}{{\rm d}x}$, so that $\vartheta=\chi(x)\delta$
with $\chi(x)=\varphi(x)/x\tfrac{{\rm d}\varphi}{{\rm d}x}(x)$.
(Compare the manipulations of Burchnall~\cite{Burchnall48}.)  Then,
the differential operators $T[e]\defeq 1+\nobreak e^{-1}\vartheta$ and
$\tilde T[\tilde e]\defeq 1+\nobreak \tilde e^{-1}\delta$ will
increment the upper parameters of the hypergeometric series in $t$
and~$x$ which define $F$ and~$\tilde F$.  That is, if one of the upper
parameters of~$F$ is~$e$, in~$T[e]F$ it will be replaced by
$1+\nobreak e$, and if none of the upper parameter is~$e$,
$T[e]F$~will have an extra parameter-pair
$\left[\begin{smallmatrix}1+e \\ e\end{smallmatrix}\right]$.  Thus,
$T[e]F$ is contiguous to~$F$ in a generalized sense.  The action of
$\tilde T[\tilde e]$ on~$\tilde F$ is similar.

The ${}_3F_2$ transformations in~\S\,\ref{subsec:extensions} are
treated as follows.  Their common form is
\begin{equation}
  F\left[
    \begin{array}{cc}
      \Delta(l+m;\,a), & (\alpha) \\
      & (\beta)
    \end{array}
    \biggm|
    t\,
    \right]
  =
  \left(1-\tfrac{x}{x_0}\right)^a
  \tilde F\left[
    \begin{array}{c:c|c}
      \begin{array}{cc}
        a, & (\gamma) \\
        & (\delta)
      \end{array}
      &
      \,
      Q_k(n)
      &
      \,
      x
    \end{array}
    \right],
  \label{eq:commonform}
\end{equation}
where the lifting function $t=\varphi(x)$ comes from
Lemma~\ref{lem:key}, i.e,
\begin{equation}
\label{eq:lifting}
  \varphi(x)=
  \varphi_{l,m;x_0}(x) \defeq
      \frac{(l+m)^{l+m}}{l^l\;m^m}
    \,
    \frac{(-x/x_0)^l}{(1-x/x_0)^{l+m}}.
\end{equation}
(Recall that $(l,m;x_0)$ is $(1,1;1)$, $(1,2;\tfrac14)$, and $(2,1;4)$
for the quadratic, first cubic, and second cubic identities.)  It
follows readily that $\chi(x) = \tfrac{x_0-x}{lx_0+mx}$, and after
some computation that $T,\tilde T$ are related by, e.g.,
\begin{equation}
\label{eq:gsrestated}
  T\left[\tfrac{a}{l+m}\right]
  = \left(1+\tfrac{m}{l} \, \tfrac{x}{x_0}
  \right)^{-1} \left(1-\tfrac{x}{x_0}\right)^{1+a}
  \tilde T\left[\tfrac{la}{l+m}\right]
  \left(1 - \tfrac{x}{x_0}\right)^{-a}.
\end{equation}
Equation~(\ref{eq:gsrestated}) is a symbolic restatement of a
well-known result of Gessel and Stanton~\cite[Prop.~2]{Gessel82}.
They applied what was essentially the operator
$T\left[\tfrac{a}{l+m}\right]$ to the classical ($k=0$) cases of
Theorems \ref{thm:A2}, \ref{thm:A3}, and~\ref{thm:A3p} to obtain their
companions: the classical ($k=0$) cases of Theorems \ref{thm:C2},
\ref{thm:C3}, and~\ref{thm:C3p}.

Now, consider the following two alternative actions on any of the
three transformation identities in~\S\,\ref{subsec:extensions} of the
form~(\ref{eq:commonform}), when $k\ge1$:
\begin{enumerate}
\item[(I)] act with $T\left[\tfrac{a}{l+m}\right]$ on it, rewriting the
  right side with the aid of~(\ref{eq:gsrestated}); or,
\item[(II)] increment $a$ (and also $b,c$ in the quadratic case), and
  decrement~$k$.
\end{enumerate}
It is easy to see that the left-hand sides resulting from actions
(I),(II) are the same, thus the resulting right sides must also be
equal.  This implies that
\begin{equation}
  \begin{split}
    \label{eq:novelty}
    &\tilde F\left[
      \begin{array}{c:c|c}
        \begin{array}{cc}
          a, & (\gamma) \\
          & (\delta)
        \end{array}
        &
        \,
        {\displaystyle \frac{n+[la/(l+m)]}{[la/(l+m)]}}
        \,
        Q_k(n)
        &
        \,
        x
      \end{array}
      \right]\\[\jot]
    &\qquad=
    \left(1+\frac{m}l \, \frac{x}{x_0}\right)
    \,
    \tilde F
    \left[
      \begin{array}{c:c|c}
        \begin{array}{cc}
          1+a, & (\gamma_+) \\
          & (\delta_+)
        \end{array}
        &
        \,
        Q_{k-1,+}(n)
        &
        \,
        x
      \end{array}
      \right],
  \end{split}
\end{equation}
where the subscript~$+$ indicates the incrementing of~$a$ (and $b,c$
in the quadratic case); and for the arrays $(\gamma)$ and~$(\delta)$,
the decrementing of~$k$ as~well.  (One sees at a glance that in all
three transformations, $(\gamma_+)=1+\nobreak (\gamma)$ and
$(\delta_+)=\nobreak (\delta)$.)  It follows by equating the
coefficients of~$x^n$ on the two sides of~(\ref{eq:novelty}) that
\begin{equation}
\label{eq:master}
\begin{split}
  K\cdot\Bigl\{{\textstyle\prod}\left(a,(\gamma)\right)\Bigr\}\, Q_k(n)
  &= A_0 \cdot\Bigl\{{\textstyle\prod}\left[n+\left(a,(\gamma)\right)\right]\Bigr\}\, Q_{k-1,+}(n)\\
  &\quad{}+ A_1\cdot \Bigl\{{\textstyle\prod}\left[(n-1) + \left(1,(\delta)\right)\right]\Bigr\}\, Q_{k-1,+}(n-1),
\end{split}
\end{equation}
with the coefficients
\begin{equation}
\label{eq:mastercoeffs}
  K=\frac{n+[la/(l+m)]}{[la/(l+m)]}, \quad A_0=1, \quad A_1=\frac{m}{l}\,\frac{1}{x_0}.
\end{equation}
Equation~(\ref{eq:master}), with~(\ref{eq:mastercoeffs}), is a
\emph{master}\/ $k$\nobreakdash-\emph{raising relation} for~$Q_k$,
standing for each of the polynomials $Q_k^{(2)}$, $Q_k^{(3)}$,
and~$Q_k^{(3')}$ of~\S\,\ref{subsec:extensions}.  It is based on a
backward difference operator on~$n$.  

By specializing $(l,m;x_0)$, one obtains an explicit
$k$\nobreakdash-raising relation for each of these three families.
For example, setting $(l,m;x_0)=(1,1;1)$ yields
\begin{equation}
\label{eq:lastmin}
\begin{split}
  abc\,\frac{n+(a/2)}{(a/2)}\,Q_k^{(2)}(n) &= (n+a)(n+b)(n+c)\,Q_{k-1,+}^{(2)}(n)\\
  &\quad{}+n(n+a-b)(n+a-c)\,Q_{k-1,+}^{(2)}(n-1)
\end{split}
\end{equation}
as the recurrence satisfied by $Q_k^{(2)}=Q_k^{(2)}(n;a;b,c)$.  This
is essentially the degree-raising relation for the dual Hahn
polynomials \cite[(9.6.8)]{Koekoek2010}.  A corresponding
$k$\nobreakdash-lowering relation can be deduced from the
hypergeometric representation of~$Q_k^{(2)}$ given in
Theorem~\ref{thm:A2}.  It is
\begin{equation}
  k\,Q_{k-1,+}^{(2)}(n) = bc
    \left(\frac{\Delta_n}{\Delta_n\lambda}\right)Q_k^{(2)}(n),
\end{equation}
where $\Delta_n$ is the forward difference operator, i.e.,
$\Delta_nf(n) = f(n+1)-f(n)$, and $\lambda\defeq n(n+a)$, so that
$\Delta_n\lambda = \allowbreak 2n+\nobreak a+\nobreak1$.  This is
equivalent to the degree-lowering relation for the dual Hahn
polynomials~\cite[(9.6.7)]{Koekoek2010}.

However, the $k$-raising relations for $Q_k^{(3)},Q_k^{(3')}$ are of a
less familiar type.  Their coefficients depend on~$k$,
unlike~(\ref{eq:lastmin}), and they do not closely resemble the
degree-raising relations for the known families of orthogonal
polynomials of a discrete argument.  This is perhaps unsurprising, as
the hypergeometric representations of $Q_k^{(3)},Q_k^{(3')}$ given in
Theorems \ref{thm:A3} and~\ref{thm:A3p} are rather novel.

The $d$-dependent ${}_4F_3$ transformations in~\S\,\ref{subsec:gens}
can be treated similarly to the ${}_3F_2$ ones
in~\S\,\ref{subsec:extensions} if the action~(II) is extended to
include an application of~$T[d]$.  In the resulting master
$k$\nobreakdash-raising relation for the $Q_k$
of~\S\,\ref{subsec:gens}, the coefficients~(\ref{eq:mastercoeffs}) are
replaced by
\begin{equation}
\label{eq:mastercoeffs2}
  K=\frac{n+[la/(l+m)]}{[la/(l+m)]}, \quad A_0=\frac{n+ld}{ld}, \quad A_1=-\,\frac{n+(a-md)}{lx_0d},
\end{equation}
which tend to the values shown in~(\ref{eq:mastercoeffs}) as
$d\to\infty$.  Setting $(l,m;x_0)=(1,1;1)$ yields
\begin{equation}
\begin{split}
  abcd\,\frac{n+(a/2)}{(a/2)}\,Q_k^{(2)}(n) &= (n+a)(n+b)(n+c)(n+d)\,Q_{k-1,+}^{(2)}(n)\\
  &\quad{}-n(n+a-b)(n+a-c)(n+a-d)\,Q_{k-1,+}^{(2)}(n-1)
\end{split}
\end{equation}
as the recurrence satisfied by the four-parameter
$Q_k^{(2)}=Q_k^{(2)}(n;a;b,c,d)$.  This is essentially the
degree-raising relation for the Racah polynomials (see
\cite[(9.2.8)]{Koekoek2010}; cf.~\cite[(3.7.6)]{Andrews99}).  A
$k$\nobreakdash-lowering formula for the four-parameter $Q_k^{(2)}$
can be deduced from its hypergeometric representation, given in
Theorem~\ref{thm:B2}, and is equivalent to the degree-lowering
relation for the Racah polynomials~\cite[(9.2.7)]{Koekoek2010}.  But
as before, the recurrences satisfied by $Q_k^{(3)},Q_k^{(3')}$ are of
a less familiar type.

The companion transformations in \S\,\ref{subsec:extcompanions} can be
treated in much the same way as those of \S\,\ref{subsec:extensions},
\emph{mutatis mutandis}.  Their common form is
\begin{equation}
  \begin{split}
    &F\left[
      \begin{array}{cc}
        \Delta(l+m;\,1+2k+a), & (\alpha)\\
        & (\beta)
      \end{array}
      \biggm|
      t\,
      \right]
    \\
    &\qquad=\left(
    1+\frac{m}{l}\,\frac{x}{x_0}
    \right)^{-1-2k}
    \left(
    1-\frac{x}{x_0}
    \right)^{1+2k+a}
    \tilde F
    \left[
      \begin{array}{c:c|c}
        \begin{array}{cc}
          a, & (\gamma) \\
          & (\delta)
        \end{array}
        &
        \,
          {\mathbf{Q}_k}(n)
          &
          \,
          x
      \end{array}
      \right].
  \end{split}
\end{equation}
To treat this form, $T\left[\frac{a}{l+m}\right]$ must be replaced
in~(I) by $T\left[\frac{1+2k+a}{l+m}\right]$, the effect of which on
each right-hand side can be worked~out by expressing it in~terms not
of~$\vartheta$ but of~$\delta$.  Also, (II)~must be replaced by its
\emph{inverse}, which decrements~$a$, etc., and increments~$k$.  By
equating the coefficients of~$x^n$ in the right-hand sides coming from
(I) and~(II), one finds after much algebraic labor an identity
resembling~(\ref{eq:master}), but with $Q_k$ replaced
by~${\mathbf{Q}}_k$ and with the new coefficient values
\begin{equation}
  K=l, \quad A_0=\frac{(l+m)n + 2lk+la}{2k+a}, \quad
  A_1=\left(\frac{m}{lx_0}\right)
  \frac{(l+m)n -2mk+la}{2k+a}.
\end{equation}
It is the master $k$\nobreakdash-raising relation for the polynomials
${\mathbf{Q}}_k^{(2)}$, ${\mathbf{Q}}_k^{(3)}$,
and~${\mathbf{Q}}_k^{(3')}$ of \S\,\ref{subsec:extcompanions}.  By
setting $(l,m;x_0)$ equal to $(1,1;1)$, $(1,2;\tfrac14)$, and
$(2,1;4)$, one obtains the relations in the statements of Theorems
\ref{thm:C2}, \ref{thm:C3}, and~\ref{thm:C3p}.  (The relation in the
first is phrased in~terms of $\hat{\mathbf{Q}}_k^{(2)}$ rather than
${\mathbf{Q}}_k^{(2)}$, but that is optional.)  This completes the
common proof of these theorems: each holds by induction on~$k$.

\section{Summation identities (I)}
\label{sec:summations1}

Besides being of intrinsic interest and of value in symbolic
manipulations, the function transformations
of~\S\,\ref{sec:maintheorems} yield new summation identities:
evaluations of (terminating) hypergeometric functions with integrally
separated and nonlinearly constrained parameters at fixed values of
their argument, such as $x=\nobreak1$.  These can be constructed by a
technique of Gessel and Stanton, which pairs transformations and their
companions.  The following lemma restates their
result~\cite[Thm.~2]{Gessel82}, which is a version of the residue
composition theorem.  (For the latter,
see~\cite[Thm.~1.2.2]{Goulden83}.)  As formulated, the lemma is
adapted to the lifting function $t=\varphi(x)=\varphi_{l,m;x_0}(x)$
defined in~(\ref{eq:lifting}).  In~it, ${\mathcal{C}}_{l,m}$ denotes
the prefactor $(l+\nobreak m)^{l+m}/l^lm^m$, and $[x^N]$~indicates the
extraction of the coefficient of~$x^N$.  Only the $l=1$, $m\ge1$ case
is stated here.


\begin{lemma}
\label{lem:GS}
  Suppose one has a pair of hypergeometric function transformations
  based on\/ $t=\varphi_{1,m;x_0}(x)$, of the form
  \begin{align*}
    G\left[
      \begin{array}{c}
        (A) \\
        (B)
      \end{array}
      \biggm|
      \varphi_{1,m;x_0}(x)
      \right]
    &=
    \left(1-\frac{x}{x_0}\right)^a
    \tilde G\left[
      \begin{array}{c}
        (\tilde A) \\
        (\tilde B)
      \end{array}
      \biggm|
      x\,
      \right],
    \\[\jot]
    G_c\left[
      \begin{array}{c}
        (A_c) \\
        (B_c)
      \end{array}
      \biggm|
      \varphi_{1,m;x_0}(x)
      \right]
    &=
    \left(1 + \frac{m}{x_0}\right)^{-1}
    \left(1-\frac{x}{x_0}\right)^{1+a_c}
    \tilde G_c\left[
      \begin{array}{c}
        (\tilde A_c) \\
        (\tilde B_c)
      \end{array}
      \biggm|
      x\,
      \right],
  \end{align*}
in which\/ $a,a_c$ appear as elements of the parameter arrays\/ $(\tilde
A),(\tilde A_c)$, respectively, and that\/ $N=(1-\nobreak a-\nobreak
a_c)/\allowbreak(1+\nobreak m)$ is a nonnegative integer.  Then,
$[x^N]\{\tilde G(x)\tilde G_c(x)\}$ equals\/
$(-{\mathcal{C}}_{1,m}/x_0)^N$ times\/ $[t^N]\{G(t)G_c(t)\}$.
Equivalently,
\begin{displaymath}
  \begin{split}
    &
    \left[
      \myatop{(\tilde A)}{(\tilde B)}
      \right]_N
    F\left[
      \begin{array}{ccc}
        -N, & (\tilde A_c), & 1-N-(\tilde B) \\
        & (\tilde B_c), & 1-N-(\tilde A)
      \end{array}
      \right]
    \\
    &
    \qquad
    =
    \left(
    -{\mathcal{C}}_{1,m}/{x_0}
    \right)^N
    \left[
      \myatop{(A)}{(B)}
      \right]_N
    F\left[
      \begin{array}{ccc}
        -N, & (A_c), & 1-N-(B) \\
        & (B_c), & 1-N-(A)
      \end{array}
      \right].
  \end{split}
\end{displaymath}
\end{lemma}

In~\cite{Gessel82}, this lemma is applied to the pair consisting of
Whipple's quadratic transformation of~${}_3F_2$ (the $k=0$ case of
Theorem~\ref{thm:A2}) and its companion (the $k=0$ case of
Theorem~\ref{thm:C2}), and yields Whipple's formula relating any very
well poised ${}_7F_6(1)$ to a $1$\nobreakdash-balanced ${}_4F_3(1)$.
(See~\cite[(5.2)]{Gessel82}.)  An extension is possible.  It can be
applied to the \emph{unrestricted} case ($k\ge0$) of
Theorem~\ref{thm:A2}, paired with the $k=0$ case of
Theorem~\ref{thm:C2}.  The lemma matches precisely the statements of
the theorems: one can read~off the hypergeometric functions $G,\tilde
G,G_c,\tilde G_c$, and their parameter arrays.  The arrays $(\tilde
A),(\tilde B)$ include $(1+\nobreak \xi_1,\dots,1+\nobreak \xi_{2k})$,
$(\xi_1,\dots,\xi_{2k})$, where $\xi_1,\dots,\xi_{2k}$ are the negated
roots of $Q_k^{(2)}(n;a,b,c)$.

To distinguish the first and second transformations of the pair, let
the parameters $a,b,c$ of Theorem~\ref{thm:A2} be renamed $d,e,f$.  In
the lemma, $a,a_c$ will accordingly signify $d,1+\nobreak a$, and the
condition that $N=(1-\nobreak a-\nobreak a_c)/(1+\nobreak
m)=0,1,2,\dots$ will become a condition that $d=-a-\nobreak 2N$.  The
left-hand $F(1)$ in the lemma is clearly a ${}_{7+2k}F_{6+2k}$, but by
examination, cancellation of parameters reduces the right-hand $F(1)$
from a ${}_6F_5(1)$ to a~${}_4F_3(1)$.

The lemma thus yields identity~(i) of the theorem below, in which the
substitutions $e\leftarrow d-\nobreak a-\nobreak N$, $f\leftarrow
e-\nobreak a-\nobreak N$ have been performed, to display a permutation
symmetry among the parameters.  The lemma can also be applied if one
replaces the ${}_3F_2$ transformation of Theorem~\ref{thm:A2} (the
first one of the pair) by the ${}_4F_3$ transformation of
Theorem~\ref{thm:B2}, in which $Q_k^{(2)}$ depends on four parameters
rather than three.  This leads to identity~(ii) of the theorem.

\begin{theorem}
  For all\/ $k\ge0$ and\/ $N\ge0$,
  the finite\ ${}_{7+2k}F_{6+2k}(1)$ sum
  {
    \setlength\arraycolsep{2.0pt}
    \begin{sizedisplaymath}{\small}
      {}_{7+2k}F_{6+2k}
      \left[
        \begin{array}{c:c}
          \begin{array}{ccccccc}
            a, & 1+\tfrac{a}2, & b, &c, &d, &e, &-N \\
            &\tfrac{a}2, &1+a-b, &1+a-c, &1+a-d, &1+a-e, &1+a+N
          \end{array}
          &
          \,
          R_k(n)
        \end{array}
        \right],
    \end{sizedisplaymath}
  }where\/ $R_k(n)$ denotes\/ $Q_k^{(2)}(N-n)/Q_k^{(2)}(N)$,
  equals\/ {\rm(i)} the finite sum
  \begin{sizealign}{\small}
    &
    \left[Q_k^{(2)}(N)\right]^{-1}\,
    \left[
      \begin{array}{c}
        1+a, \:\: 1-k+a-d-e \\
        1+a-d, \:\: 1+a-e
      \end{array}
      \right]_N
    \nonumber
    \\
    &
    \qquad{}\times{}_4F_3
    \left[
      \begin{array}{c}
        1+a-b-c, \:\: d, \:\: e, \:\: -N\\
        1+a-b, \:\: 1+a-c, \:\: k-a+d+e-N
      \end{array}
      \right]
    \nonumber
    \end{sizealign}
  if\/ $Q_k^{(2)}(n)\defeq Q_k^{(2)}(n;-a-2N; d-a-N, e-a-N)$, and\/ {\rm(ii)}
  the finite sum
  \begin{sizealign}{\small}
    &
    \left[Q_k^{(2)}(N)\right]^{-1}\,
    \left[
      \begin{array}{c}
        1+a, \:\: 1-k+a-d-e, \:\: 1-k+a-f \\
        1+a-d, \:\: 1+a-e, \:\: 1+a-f
      \end{array}
      \right]_N
    \nonumber
    \\
    & \qquad{}\times {}_5F_4
    \left[
      \begin{array}{c}
        1+a-b-c, \:\: d, \:\: e, \:\: 1+a-f, \:\: -N\\
        1+a-b, \:\: 1+a-c, \:\: k-a+d+e-N, \:\: 1-k+a-f
      \end{array}
      \right]
    \nonumber
  \end{sizealign}
  if\/ $Q_k^{(2)}\defeq Q_k^{(2)}(n;\allowbreak -a-2N;
  d-a-N, e-a-N, f-a-N)$.
\end{theorem}

The two identities of the theorem reduce to Whipple's formula when
$k=0$, and the second reduces to the first when $f\to\infty$.  In
both, the left-hand ${}_7F_6(1)$ is very well poised and the
right-hand series has $S=1+\nobreak k$, resp.\ $S=1$.  The $k=1$ case
of the second identity can be shown to agree with a result of
Srivastava, Vyas and Fatawat~\cite[Thm.~3.2]{Srivastava2018} by using
the formula (\ref{eq:Q1}) for the four-parameter quadratic polynomial
$Q_1^{(2)}$.  Like Whipple's formula
(cf.\ \cite[\S\S\,3.4,\,3.5]{Andrews99}), they have interesting
specializations and limits.  For example, if the
${}_{7+2k}F_{6+2k}(1)$ has $S=2$, then the right-hand parameters
$1+\nobreak a-\nobreak b-\nobreak c$, $k-\nobreak a+\nobreak
d+\nobreak e-\nobreak N$ will equal each other and can be cancelled.
The two identities then become extensions of Dougall's theorem on the
sum of a $2$\nobreakdash-balanced, very well poised ${}_7F_6(1)$~
\cite[Thm.~3.5.1]{Andrews99}.

One can also apply Lemma~\ref{lem:GS} to the pair consisting of the
unrestricted Theorem~\ref{thm:A3} (the first cubic transformation
of~${}_3F_2$), resp.\ Theorem~\ref{thm:B3} (the first cubic
transformation of~${}_4F_3$), and the $k=0$ case of its companion,
Theorem~\ref{thm:C3}.  The two summation identities which result are
extensions to $k\ge0$ of the first cubic summation identity of Gessel
and Stanton~\cite[(1.7)]{Gessel82}.  Details are left to the reader.

\section{Summation identities (II)}
\label{sec:summations2}

One can obtain a parametric finite summation identity from any of the
extended function transformations of~\S\,\ref{sec:maintheorems} by a
classical technique: multiplying both sides by a power of
$1-\nobreak{x}$, such as $(1-\nobreak x)^{w-a+m-1}$, and equating the
coefficients of~$x^m$ on the two sides.  This technique was applied by
Bailey to many hypergeometric transformations, including Whipple's
quadratic transformation of ${}_3F_2$ and its companion (the $k=0$
cases of Theorems \ref{thm:A2} and~\ref{thm:C2}); see \cite[p.~97,
  Examples~5,6]{Bailey35}.  Applying it to the unrestricted ($k\ge0$)
versions of Theorems \ref{thm:A2}, \ref{thm:B2}, and~\ref{thm:C2} is
straightforward and yields:

\begin{theorem}
\label{thm:71}
  For all\/ $k\ge0$ and\/ $m\ge0$, one has\/ {\rm(i)} the finite
  summation identity
  \begin{align*}
    &
    {}_5F_4\left[
      \begin{array}{c:c}
        \begin{array}{c}
          \tfrac{a}2, \:\: \tfrac12+\tfrac{a}2, \:\: 1-k+a-b-c, \\
          1+a-b, \:\: 1+a-c,
        \end{array}
        &
        \begin{array}{cc}
          1+a-w, & -m \\
          \tfrac12+\tfrac12(a-w-m), & 1+\tfrac12(a-w-m)
        \end{array}
      \end{array}
      \right]
    \\
    &{}\qquad=
    \left[
      \begin{array}{c}
        w\\
        w-a
      \end{array}
      \right]_m
    \:
      {}_{4+2k}F_{3+2k}
      \left[
        \begin{array}{c:c:c}
          \begin{array}{c}
            a, \:\: b, \:\: c, \\
            1+a-b, \:\: 1+a-c,
          \end{array}
          &
          \begin{array}{c}
            -m \\
            w
          \end{array}
          &
          \,
          Q_k^{(2)}(n)
        \end{array}
        \right],
  \end{align*}
where\/ $Q_k^{(2)}(n)\defeq Q_k^{(2)}(n;a;b,c)$, and\/ {\rm (ii)} a
like formula in which a parameter-pair\/
$\left[\begin{smallmatrix}k+d \\ d\end{smallmatrix}\right]$ is added
  to the left-hand side, and\/ $Q_k^{(2)}(n)\defeq
  Q_k^{(2)}(n;a;b,c,d)$.
\end{theorem}
\begin{theorem}
\label{thm:72}
  For all\/ $k\ge0$ and\/ $m\ge0$, one has the finite summation identity
  \begin{sizemultline}{\small}
    {}_{5+k}F_{4+k}\biggl[
      \begin{array}{c:}
        \begin{array}{c}
          \tfrac12+k+\tfrac{a}2, \:\: 1+k+\tfrac{a}2, \:\: 1-k+a-b-c, \\
          1+a-b, \:\: 1+a-c,
        \end{array}
      \end{array}
      \nonumber
      \\
      \begin{array}{c:c}
        \begin{array}{cc}
          1+a-w, & -m \\
          1+\tfrac12(a-w+2k-m), & \tfrac32+\tfrac12(a-w+2k-m)
        \end{array}
        &
        \,
        P_k(n)
      \end{array}
      \biggr]
    \nonumber
    \\
    \quad=
    \frac{(1+a-w)_{1+2k}(w)_m}{(w-a-1-2k)_m}
    \:\:
      {}_{5+4k}F_{4+4k}
      \left[
        \begin{array}{c:c:c}
          \begin{array}{c}
            a, \:\: b, \:\: c, \\
            1+a-b, \:\: 1+a-c,
          \end{array}
          &
          \begin{array}{c}
            -m \\
            w
          \end{array}
          &
          \,
          {\mathbf{Q}}_k^{(2)}(n)
        \end{array}
        \right],
  \end{sizemultline}
  where\/ ${\mathbf{Q}}_k^{(2)}(n)\defeq
  {\mathbf{Q}}_k^{(2)}(n;a;b,c)$, which is of degree\/ $1+4k$ in\/ $n$, and\/ $P_k(n)\defeq
  P_k(n;1+a-w,-m)$ is a polynomial of degree\/ $k$ in~\/$n$ defined by
  \begin{displaymath}
    P_k(n;\,A,B) \defeq
    (n+A)_{2k+1}\:\,
   {}_2F_1
  \left[
    \begin{array}{c|c}
      \begin{array}{c}
      {-1-2k,\:\: n+B} \\
      {-n-A-2k}
      \end{array}
      &
      -1
    \end{array}
  \right],
  \end{displaymath}
  which\/ {\rm(}by series reversal\,{\rm)} is odd under the interchange
  of\/ $A,B$.
\end{theorem}

In the left-hand ${}_{5+k}F_{4+k}(1)$ of Theorem~\ref{thm:72}, the
convention introduced in~\S\,\ref{sec:prelims} is not adhered~to, for
simplicity of expression: the weighting function, here $P_k(n) =
P_k(n;A,B)$, does not equal unity at~$n=0$.  For instance, $P_0(n)$ is
identically equal to $A-B$.  It is worth mentioning that the
polynomials $P_k(n;A,B)$ have the generating function
\begin{equation}
\begin{split}
&  \sum_{k=0}^\infty P_k(n;\,A,B)\,\frac{t^{1+2k}}{(1+2k)!} \\
& \qquad = \frac12\left[
(1-t)^{-n-A}(1+t)^{-n-B} 
-
(1-t)^{-n-B}(1+t)^{-n-A}
\right],
\end{split}
\end{equation}
which is a specialization of the Srivastava--Singhal generating
function for Jacobi polynomials~\cite{Srivastava73}.

The summation identities in Theorems \ref{thm:71} and~\ref{thm:72}
reduce when $k=0$ to those given by Bailey
\cite[\S\,4.5(1,2)]{Bailey35}.  In each, the left-hand series
(${}_5F_4(1)$ or ${}_6F_5(1)$, resp.\ ${}_{5+k}F_{4+k}(1)$) is either
$(1+\nobreak k)$-balanced or 1\nobreakdash-balanced, and the
right-hand one (${}_{4+2k}F_{3+2k}(1)$, resp.\ ${}_{5+4k}F_{4+4k}(1)$)
is nearly poised.  The $k=1$ case of Theorem~\ref{thm:71}(ii) was
recently proved by Wang and Rathie~\cite[Cor.~4]{Wang2013}.  There are
some interesting specializations of the $k\ge1$ cases of Theorems
\ref{thm:71} and~\ref{thm:72}, which can be viewed as extensions of
Bailey's several specializations of the $k=0$ case of
Theorem~\ref{thm:71}(i) (for the latter,
see~\cite[\S\,4.5]{Bailey35}).

Bailey noted that there is an equivalence between Whipple's quadratic
transformation of~${}_3F_2$, i.e., the $k=0$ case of
Theorem~\ref{thm:A2}, and his formula relating a
1\nobreakdash-balanced ${}_5F_4(1)$ to a nearly poised ${}_4F_3(1)$,
i.e., the $k=0$ case of Theorem~\ref{thm:71}(i): one implies the
other.  (For a $q$\nobreakdash-analogue, see~\cite{AlSalam84}.)  One
now sees that this equivalence holds in greater generality, in a
manner parametrized by $k=0,1,2,\dots$.

\section*{Acknowledgements}

We thank the anonymous referees for several corrections and
enhancements.


\begin{thebibliography}{10}
\expandafter\ifx\csname url\endcsname\relax
  \def\url#1{\texttt{#1}}\fi
\expandafter\ifx\csname urlprefix\endcsname\relax\def\urlprefix{URL }\fi
\expandafter\ifx\csname href\endcsname\relax
  \def\href#1#2{#2} \def\path#1{#1}\fi

\bibitem{AlSalam84}
W.~A. Al-Salam, A.~Verma, On quadratic transformations of basic series, {SIAM}
  J.~Math. Anal. 15~(2) (1984) 414--421.

\bibitem{Andrews99}
G.~E. Andrews, R.~Askey, R.~Roy, Special Functions, Vol.~71 of The Encyclopedia
  of Mathematics and its Applications, Cambridge Univ. Press, Cambridge, UK,
  1999.

\bibitem{Askey94}
R.~Askey, A look at the {B}ateman project, in: W.~Abikoff, J.~S. Birman,
  K.~Kuiken (Eds.), The Mathematical Legacy of {W}ilhelm {M}agnus: Groups,
  Geometry, and Special Functions, Vol. 169 of Contemporary Mathematics,
  American Mathematical Society, Providence, RI, 1994, pp. 29--43.

\bibitem{Bailey28}
W.~N. Bailey, Products of generalized hypergeometric series, Proc. London Math.
  Soc.~(2) 28~(4) (1928) 242--254.

\bibitem{Bailey35}
W.~N. Bailey, Generalized Hypergeometric Series, no.~32 in Cambridge Tracts in
  Mathematics and Mathematical Physics, Cambridge Univ. Press, Cambridge, UK,
  1935.

\bibitem{Burchnall48}
J.~L. Burchnall, On the well-poised {${}_3F_2$}, J.~London Math. Soc. 23 (1948)
  253--257.

\bibitem{Chen2005}
Kung-Yu Chen, H.~M. Srivastava, Series identities and associated families of
  generating functions, J.~Math. Anal. Appl. 311~(2) (2005) 582--599.

\bibitem{Gessel82}
I.~Gessel, D.~Stanton, Strange evaluations of hypergeometric series, {SIAM}
  J.~Math. Anal. 13~(2) (1982) 295--308.

\bibitem{Goulden83}
I.~P. Goulden, D.~M. Jackson, Combinatorial Enumeration, Wiley, New York, 1983.

\bibitem{Goursat1881}
{\'E}.~Goursat, Sur l'{\'e}quation diff{\'e}rentielle lin{\'e}aire qui admet
  pour int{\'e}grale la s{\'e}rie hyperg{\'e}om{\'e}trique, Ann. Sci. {\'Ec}.
  Norm. Sup.~(2) 10 (1881) S3--S142.

\bibitem{Karlsson71}
P.~W. Karlsson, Hypergeometric functions with integral parameter differences,
  J.~Math. Phys. 12 (1971) 270--271.

\bibitem{Kato2008}
M.~Kato, Algebraic transformations of {${}_3F_2$}, Funkcial. Ekvac. 51~(2)
  (2008) 221--243.

\bibitem{Koekoek2010}
R.~Koekoek, P.~A. Lesky, R.~F. Swarttouw, Hypergeometric Orthogonal Polynomials
  and Their {$q$}-Analogues, Springer-Verlag, Berlin, 2010, with a foreword by
  T.~H. Koornwinder.

\bibitem{Miller2013}
A.~R. Miller, R.~B. Paris, Transformation formulas for the generalized
  hypergeometric function with integral parameter differences, Rocky Mountain
  J.~Math. 43~(1) (2013) 291--327.

\bibitem{Niblett51}
J.~D. Niblett, Some hypergeometric identities, Pacific J.~Math. 2~(2) (1953)
  219--225.

\bibitem{Rakha2011}
M.~A. Rakha, A.~K. Rathie, Extensions of {E}uler type {II} transformations and
  {S}aalsch{\"u}tz's theorem, Bull. Korean Math. Soc. 48~(1) (2011) 151--156.

\bibitem{Rakha2009}
M.~A. Rakha, N.~Rathie, P.~Chopra, On an extension of a quadratic
  transformation formula due to {K}ummer, Math. Commun. 14~(2) (2009) 207--209.

\bibitem{Srivastava84}
H.~M. Srivastava, H.~L. Manocha, A Treatise on Generating Functions, Halsted
  Press, New York, 1984.

\bibitem{Srivastava73}
H.~M. Srivastava, J.~P. Singhal, New generating functions for {J}acobi and
  related polynomials, J.~Math. Anal. Appl. 41~(3) (1973) 748--752.

\bibitem{Srivastava2018}
H.~M. Srivastava, Y.~Vyas, K.~Fatawat, Extensions of the classical theorems for
  very well-poised hypergeometric functions, Rev. R.~Acad. Cienc. Exactas Fís.
  Nat. Ser.~A Math. RACSAM (2019) in press, previous version available as
  {\path{arXiv:1410.3241}}, current version available as doi:10.1007/s13398-017-0485-5.

\bibitem{Vidunas2009}
R.~Vid{\=u}nas, Algebraic transformations of {G}auss hypergeometric functions,
  Funkcial. Ekvac. 52 (2009) 139--180,
  {\path{arXiv:math.CA/0408269}}.

\bibitem{Wang2013}
Xiaoxia Wang, A.~K. Rathie, Extension of a quadratic transformation due to
  {W}hipple, with an application, Adv. Difference Equ. 103
  (2013) article id~157 [8~pp.].

\end{thebibliography}

\end{document}